\pdfoutput=1
\RequirePackage{ifpdf}
\ifpdf 
\documentclass[pdftex]{sigma}
\else
\documentclass{sigma}
\fi

\usepackage{enumitem}

\numberwithin{equation}{section}

\newtheorem{Theorem}{Theorem}[section]
\newtheorem*{Theorem*}{Theorem}
\newtheorem{Corollary}[Theorem]{Corollary}
\newtheorem{Lemma}[Theorem]{Lemma}
\newtheorem{Proposition}[Theorem]{Proposition}
 { \theoremstyle{definition}
\newtheorem{Definition}[Theorem]{Definition}

\newtheorem{Example}[Theorem]{Example}
\newtheorem{Remark}[Theorem]{Remark} }

\usepackage{tikz-cd}

\usepackage[all]{xy}

\newcommand{\cat}[1]{\mathcal{#1}}
\newcommand{\SV}{{\operatorname{SVec}}}

\newcommand{\FPdim}{\operatorname{FPdim}}
\def\cross#1#2{{#1}_{#2}^{\times}}
\newcommand{\cb}{\cat{B}}

\newcommand{\Irr}{\operatorname{Irr}}
\newcommand\Vc{\operatorname{Vec}}

\newcommand{\unit}{\mathbf{1}}

\newcommand{\id}{{\operatorname{id}}}

\newcommand{\Id}{{\operatorname{Id}}}
\newcommand{\Rep}{{\operatorname{Rep}}}

\newcommand{\Autr}[1]{\operatorname{Aut}_{\otimes}(#1)}

\newcommand{\C}{\mathcal{C}}
\def\Pic#1{\text{Pic}(#1)}

\def\uu#1{\underline{\underline{#1}}}
\newcommand{\vg}{\vec{g}}
\newcommand{\va}{\vec{a}}
\newcommand{\vb}{\vec{b}}

\def\cB{\mathcal{B}}
\def\cE{\mathcal{E}}

\def\cL{\mathcal{L}}
\def\cN{\mathcal{N}}
\def\cM{\mathcal{M}}

\def\cC{\mathcal{C}}

\def\cS{\mathcal{S}}
\def\cD{\mathcal{D}}

\begin{document}

\allowdisplaybreaks

\newcommand{\arXivNumber}{1712.07097}

\renewcommand{\PaperNumber}{085}

\FirstPageHeading

\ShortArticleName{Categorical Fermionic Actions and Minimal Modular Extensions}

\ArticleName{Categorical Fermionic Actions and Minimal Modular\\ Extensions}

\Author{C\'{e}sar GALINDO and C\'esar F. VENEGAS-RAM\'IREZ}

\AuthorNameForHeading{C.~Galindo and C.~Venegas-Ram\'irez}

\Address{Departamento de Matem\'aticas, Universidad de los Andes, Bogot\'a, Colombia}
\Email{\href{cn.galindo1116@uniandes.edu.co}{cn.galindo1116@uniandes.edu.co}, \href{cf.venegas10@uniandes.edu.co}{cf.venegas10@uniandes.edu.co}}

\ArticleDates{Received March 12, 2025, in final form October 04, 2025; Published online October 13, 2025}

\Abstract{We define fermionic actions of finite super-groups on fermionic fusion categories and establish necessary and sufficient conditions for their existence. Our main result characterizes when a braided fusion category admits a minimal non-degenerate extension in terms of cohomological obstructions. This characterization for braided fusion categories with non-Tannakian M\"uger center involves the fermionic structures and fermionic actions introduced in this work.}


\Classification{18M20}

\section{Introduction}

The purpose of this paper is to study minimal non-degenerate extensions of braided fusion categories, with an emphasis on minimal non-degenerate extensions of non-Tannakian symmetric fusion categories.

Let $\cB$ be a premodular fusion category. Minimal modular extensions were defined by Michael M\"uger in \cite{Mu2} as a modular category $\cM \supset \cB$ such that the centralizer of $\mathcal{Z}_2(\cB)$ in $\cM$ is $\cB$, where $\mathcal{Z}_2(\cB)$ is the M\"uger center of $\cB$. We work with the natural generalization of minimal non-degenerate extensions, where we only require $\cM$ to be non-degenerate (see Definition~\ref{def: minimal non-degenerate extension}). This broader notion coincides with M\"uger's definition in the unitary case but provides additional flexibility for analyzing non-unitary braided fusion categories.

Braided fusion categories play a fundamental role in diverse areas of mathematics and theoretical physics. They appear in the representation theory of quantum groups and Hopf algebras~\cite{CP94}, the construction of invariants of 3-manifolds \cite{RT}, and as a mathematical framework for topological phases of matter \cite{MS89, Wen90, Wit89}. In the context of condensed matter physics, braided fusion categories organize the particle excitations in (2+1)-dimensional topological quantum systems, leading to their importance in quantum computation \cite{Kit03, Wan10}.

Following \cite{LKW}, we denote the set of equivalence classes of minimal non-degenerate extensions of $\cB$ by $\cM_{\operatorname{ext}}(\cB)$. In \cite{LKW}, it was proved that if $\cE$ is a symmetric fusion category, $\cM_{\operatorname{ext}}(\cE)$ has a natural abelian group structure. For a general premodular category, $\cM_{\operatorname{ext}}(\cB)$ is a torsor (possibly empty) over $\cM_{\operatorname{ext}}(\mathcal{Z}_2(\cB))$. Hence, the results in \cite{LKW} reduce the problem of constructing all minimal non-degenerate extensions of a premodular category $\cB$ to the following steps:
\begin{enumerate}\itemsep=0pt
 \item[(1)]\label{primer paso} Determine if $\cM_{\operatorname{ext}}(\cB)$ is empty or not. If $\cM_{\operatorname{ext}}(\cB)$ is not empty, construct one minimal extension of $\cB$.
 \item[(2)]\label{segundo paso} Construct all minimal extensions of the symmetric fusion category $\mathcal{Z}_2(\cB)$.
\end{enumerate}

The goal of this paper is to develop an obstruction theory for the existence of minimal non-degenerate extensions of braided fusion categories, with particular emphasis on the non-Tannakian case involving super-groups.

The main results of the present paper are the following.

For braided fusion categories $\mathcal{B}$ whose M\"uger center is Tannakian (i.e., equivalent to $\Rep(G)$ for some finite group $G$), we show that the obstruction to having a minimal non-degenerate extension is an element $O_4(\mathcal{B})$ in $H^4(G; \mathbb{C}^\times)$, which we call the $H^4$-anomaly of $\mathcal{B}$. This obstruction has a concrete formula in many cases (see Corollary~\ref{formula 4 quasi-trivial}), allowing for explicit calculations. Using this formula, we provide specific examples of braided fusion categories without minimal non-degenerate extensions in Section \ref{subsection:obstruction 4}, including a noteworthy example first discovered by Drinfeld.

For braided fusion categories $\mathcal{B}$ whose M\"uger center is super-Tannakian \big(i.e., equivalent to~\smash{$\Rep\bigl(\widetilde{G},z\bigr)$} for some finite super-group \smash{$\bigl(\widetilde{G},z\bigr)$\big)}, we establish a framework involving fermionic actions of super-groups. We introduce the concept of categorical fermionic actions of super-groups on fermionic fusion categories (see Definition~\ref{def:fermionic-action}). This generalization provides the necessary tools to analyze minimal non-degenerate extensions in the non-Tannakian case.

To state our main result precisely, we introduce the key constructions involved. Given a~braided fusion category $\cB$ with M\"uger center $\mathcal{Z}_2(\cB) = \smash{\Rep\bigl(\widetilde{G},z\bigr)}$ for a super-group \smash{$\bigl(\widetilde{G},z\bigr)$}, we denote by \smash{$G := \widetilde{G}/\langle z \rangle$} the quotient by the central element $z$ of order two. The maximal central Tannakian subcategory of $\cB$ is then equivalent to $\Rep(G)$. We denote by $\cB_G$ the \emph{de-equivariantization} of $\cB$ with respect to $\Rep(G)$ (see Section~\ref{categorical actions}), which carries a canonical spin-braided structure with fermion given by the generator of $\mathcal{Z}_2(\cB_G) \cong \SV$ (see \cite[Proposition~4.30\,(iii)]{DGNO}).

A \emph{fermionic action} of \smash{$\bigl(\widetilde{G},z\bigr)$} on a spin-braided fusion category consists of a categorical $G$-action together with additional compatibility data encoding the super-group structure (see Definition~\ref{def:fermionic-action} for the precise definition). Crucially, a fermionic $\bigl(\widetilde{G},z\bigr)$-action on a category $\cS$ induces an underlying categorical $G$-action on~$\cS$.

Our main result, which provides a characterization of when minimal non-degenerate extensions exist, is the following.

\begin{Theorem}[Theorem~\ref{theorem:teorema de-equivariantización y extensiones}]\label{thm:main-intro}
Let $\cB$ be a braided fusion category with non-trivial maximal central Tannakian subcategory $\Rep(G) \subseteq \mathcal{Z}_2(\cB)$.
\begin{enumerate}\itemsep=0pt
 \item[$(1)$] If $\mathcal{Z}_2(\cB) = \Rep(G)$ $($the \emph{modularizable} case$)$, then $\cB$ admits a minimal non-degenerate extension if and only if the $H^4$-anomaly $O_4(\cB) \in H^4(G,\mathbb{C}^\times)$ vanishes.
 \item[$(2)$] If $\mathcal{Z}_2(\cB) = \Rep\bigl(\widetilde{G},z\bigr)$ for a super-group $\bigl(\widetilde{G},z\bigr)$ with $G = \widetilde{G}/\langle z \rangle$ $($the \emph{non-modularizable} case$)$, then $\cB$ admits a minimal non-degenerate extension if and only if the following three conditions hold:
 \begin{enumerate}\itemsep=0pt
 \item[$(a)$] The de-equivariantization $\cB_G$ $($which is slightly degenerate with $\mathcal{Z}_2(\cB_G) \cong \SV)$ has a minimal non-degenerate extension $\cS$.

 \item[$(b)$] There exists a fermionic action of $\bigl(\widetilde{G},z\bigr)$ on $\cS$ whose underlying $G$-action stabilizes~${\cB_G \subseteq \cS}$ and restricts to the canonical $G$-action on $\cB_G$ induced by de-equivarianti\-za\-tion.

 \item[$(c)$] The $H^4$-anomaly $O_4\bigl(\cS^G\bigr)$ of the $G$-equivariantization $\cS^G$ vanishes.
 \end{enumerate}
\end{enumerate}

Here the $H^4$-anomaly is the cohomological obstruction defined in Definition~{\rm\ref{def: H4 anolmaly}}, and fermionic actions of super-groups are introduced in Section {\rm\ref{acciosuper}}.
\end{Theorem}

\begin{Remark}
Part (1) provides an obstruction theory for the modularizable case, reducing the problem to the vanishing of a single cohomology class. Part (2) reveals the structure of the non-modularizable case: one must first construct a minimal extension $\cS$ of the slightly degenerate category $\cB_G$, then verify that $\cS$ admits a fermionic super-group action compatible with the original $G$-action on $\cB_G$. The necessity of fermionic structures in part (2) motivates our systematic development of fermionic action theory in Section \ref{acciosuper}.

The condition in part (2)(b) is subtle: not every $G$-action on $\cS$ extends to a fermionic $\bigl(\widetilde{G},z\bigr)$-action, and even when such extensions exist, they may not preserve the subcategory $\cB_G$. This compatibility requirement is encoded in the fermionic action structure and will be analyzed in detail using cohomological obstructions in Theorem~\ref{theorem:bstruction fermionicaction}.
\end{Remark}

Johnson-Freyd and Reutter \cite{JFR} recently demonstrated that every slightly degenerate braided fusion category admits a minimal nondegenerate extension, thereby resolving a longstanding open problem. Their analysis, based on fusion 2-categories, shows that the relevant obstruction—lying in $H^5\bigl(K(\mathbb{Z}_2,2);\mathbb{C}^\times\bigr)\cong \mathbb{Z}_2$—vanishes universally. This result, combined with our obstruction theory for the Tannakian case, significantly advances the classification program for minimal nondegenerate extensions of braided fusion categories. Our work complements these findings by clarifying the conditions for existence and elucidating the nature of the obstructions, offering a perspective that integrates with the result of \cite{JFR}.

The paper is organized as follows:

In Section \ref{section: preliminaries}, we recall necessary definitions and results about fusion categories, braided and symmetric fusion categories, and actions of groups on fusion categories.

Section \ref{acciosuper} develops the theory of fermionic actions on fermionic fusion categories and spin-braided fusion categories. In Theorem~\ref{theorem: equivalencia dequivariantizacion- equivariantizacion} and Corollary~\ref{coro: equivalencia dequivariantizacion- equivariantizacion braided case}, we prove a 2-equivalence between the 2-category of fermionic fusion categories with fermionic actions of a super-group \smash{$\bigl(\widetilde{G},z\bigr)$} and the 2-category of fusion categories over \smash{$\Rep\bigl(\widetilde{G},z\bigr)$}, analogous to \cite[Theorem~4.18]{DGNO}. We define a~cohomological obstruction to the existence of fermionic actions and provide a group-theoretical interpretation of this obstruction. Finally, we present results about fermionic actions on non-degenerate spin-braided fusion categories of dimension four.

In Section \ref{section: minimal non-degenerated extension}, we study minimal non-degenerate extensions of braided fusion categories, with emphasis on minimal non-degenerate extensions of super-Tannakian fusion categories. In Section~\ref{subsection: Obstruction theory to existence of minimal non-degenerate extensions}, we establish necessary and sufficient conditions for the existence of a minimal non-degenerate extension of a fusion category. Additionally, in Corollary~\ref{corol:Ostrik}, we prove that the group homomorphism between the minimal non-degenerate extensions of $\Rep\bigl(\widetilde{G},z\bigr)$ and $\SV$ is surjective if and only if the super-group \smash{$\bigl(\widetilde{G},z\bigr)$} is trivial. Finally, in Section~\ref{subsection:obstruction 4}, we present examples of braided fusion categories without minimal non-degenerate extensions. Continuing with the study of fermionic actions on pointed spin-braided fusion categories of rank four, we investigate their $H^4$ obstruction.

\section{Preliminaries} \label{section: preliminaries}

We will briefly review the basic definitions and results of fusion categories and braided fusion categories that are necessary for defining fermionic actions and study minimal non-degenerate extensions. For comprehensive treatments, we refer to \cite{DGNO,Book-ENO}.

\subsection{Fusion categories}

A \emph{fusion category} (over $\mathbb{C}$) is a $\mathbb{C}$-linear semisimple rigid tensor category $\cC$ with finitely many isomorphism classes of simple objects, finite-dimensional spaces of morphisms, and such that the unit object $\mathbf{1}$ is simple. All fusion categories considered in this paper are over $\mathbb{C}$. By a \emph{fusion subcategory} of a fusion category, we mean a full tensor abelian subcategory.

A fusion category is called \emph{pointed} if all its simple objects are invertible with respect to the tensor product. For a fusion category $\cC$, we denote by $\cC_{\rm pt}$ the maximal pointed fusion subcategory of $\cC$, consisting of all invertible objects and their direct sums.

Up to equivalence, every pointed fusion category has the form $\operatorname{Vec}_G^\omega$, where $G$ is a finite group and $\omega \in Z^3(G,\mathbb{C}^\times)$ is a 3-cocycle. This is the category of $G$-graded vector spaces with associator twisted by $\omega$. The simple objects are denoted by $\mathbb{C}_g$ for $g \in G$, with tensor product~${\mathbb{C}_g \otimes \mathbb{C}_h = \mathbb{C}_{gh}}$. The associativity constraint is given by
\[
a_{g,h,k} = \omega(g,h,k) \cdot \id_{\mathbb{C}_{ghk}} \colon\ (\mathbb{C}_g \otimes \mathbb{C}_h) \otimes \mathbb{C}_k \to \mathbb{C}_g \otimes (\mathbb{C}_h \otimes \mathbb{C}_k)
\]
for all $g,h,k \in G$. The 3-cocycle condition ensures that the associativity constraint satisfies the pentagon axiom.

The group of isomorphism classes of invertible objects in $\cC$ is denoted by $\operatorname{Inv}(\cC)$.

We denote by $\Irr(\cC)$ the set of isomorphism classes of simple objects in $\cC$ and by $K_0(\cC)$ the \emph{Grothendieck ring} of $\cC$. The \emph{rank} of $\cC$ is $|{\Irr(\cC)}|$.

There exists a unique ring homomorphism $\FPdim\colon K_0(\cC) \to \mathbb{R}$ (see \cite[Proposition~3.3.6]{Book-ENO}) such that $\FPdim(X) > 0$ for any $X \in \Irr(\cC)$. The \emph{Frobenius--Perron dimension} of a fusion category $\cC$ is defined as
$\FPdim(\cC) = \sum_{X \in \Irr(\cC)} \FPdim(X)^2$.

Let $F \colon \cC \to \cD$ be a monoidal functor between fusion categories. The \emph{image} of $F$, denoted~$\operatorname{Im}(F)$, is the smallest fusion subcategory of $\cD$ containing all objects of the form $F(X)$ for $X \in \cC$. If $F$ is faithful (that is, the induced functor $\cC \to \operatorname{Im}(F)$ is an equivalence), then $\FPdim(\cC) \leq \FPdim(\cD)$ with equality if and only if $F$ is an equivalence \cite[Proposition~6.3.3]{Book-ENO}. Conversely, if $F$ is essentially surjective (meaning $\operatorname{Im}(F) = \cD$), then $\FPdim(\cC) \geq \FPdim(\cD)$ with equality if and only if $F$ is an equivalence \cite[Proposition~6.3.4]{Book-ENO}.
\subsection{Braided fusion categories}\label{section: braided fusion categories}

A fusion category $\cB$ is called \emph{braided} if it is endowed with a natural isomorphism
\begin{align*}
c_{X,Y} \colon\ X \otimes Y\to Y \otimes X, \qquad X,Y \in\cB,
\end{align*}
satisfying the hexagon axioms, see \cite{JS}. These axioms ensure that the braiding is compatible with the tensor structure of the category.

If $\cD$ is a full subcategory of $\cB$, the \emph{centralizer} of $\cD$ with respect to $\cB$ is defined as the full subcategory
\[ C_{\cB}(\cD):=\{ Y \in \cB \mid c_{Y,X} \circ c_{X,Y}=\id_{X \otimes Y} \text{ for all } X \in \cD \}.\]

The \emph{M\"uger center} of $\cB$ is the fusion subcategory $\mathcal{Z}_2(\cB):= C_{\cB}(\cB)$, that is,
\[\mathcal{Z}_2(\cB)=\{Y \in \cB\mid c_{Y,X} \circ c_{X,Y}=\id_{X \otimes Y }, \ \text{for all } X \in \cB\}.\]

A braided fusion category $\cb$ is called \emph{symmetric} if $\mathcal{Z}_2(\cb)=\cb$, i.e., if $c_{Y,X} \circ c_{X,Y}=\id_{X \otimes Y}$ for each pair of objects $X,Y$ in $\cb$.

Symmetric fusion categories are equivalent to one of the following two examples:
\begin{enumerate}[leftmargin=*,label={\rm (\alph*)}]\itemsep=0pt
 \item {\it Tannakian categories}. The category $\operatorname{Rep}(G)$ of finite dimensional complex representation of a finite group $G$, with standard braiding $c_{X,Y}(x \otimes y):= y \otimes x$ for $x \in X$ and $y \in Y$.

 \item {\it Super-Tannakian categories}. A \emph{finite super-group} is a pair $(G,z)$, where $G$ is a finite group and $z$ is a central element of order two. An irreducible representation of $G$ is called \emph{odd} if $z$ acts as the scalar $-1$, and \emph{even} if $z$ acts as the identity. The \emph{degree} of an irreducible representation is $0$ if it is even, and $1$ if it is odd. If the degree of a simple object $X$ is denoted by $|X| \in \{0,1\}$, then the braiding of two simple objects $X$, $Y$ is
\[ c'_{X,Y}(x\otimes y)= (-1)^{|X||Y|}y\otimes x.\]
The category $\operatorname{Rep}(G)$ with the braiding $c'$ is called a \emph{super-Tannakian} category, and will be denoted by $\Rep(G,z)$.
\end{enumerate}

The super-Tannakian category $\operatorname{Rep}(\mathbb{Z}/2\mathbb{Z},1)$ is called the category of \emph{super-vector spaces} and will be denoted by $\SV$. This category will play a central role in our study of fermionic actions.

In \cite{deligne2002categories}, Deligne establishes that every symmetric fusion category is braided equivalent to $\operatorname{Rep}(G)$ or $\operatorname{Rep}(G,z)$ for a unique finite group $G$ or super group $(G,z)$.

A braided fusion category $(\cB,c)$ is called \emph{non-degenerate} if $\mathcal{Z}_2(\cB) \cong \Vc$, that is, its M\"uger center is trivial. For \emph{spherical} braided fusion categories, non-degeneracy is equivalent to modularity, meaning the invertibility of the $S$-matrix.

\begin{Example}[pointed braided fusion categories]
Let $A$ be a finite abelian group. A braided structure on $\operatorname{Vec}_A^{\omega}$ is determined by a function $c\colon A \times A \to \mathbb{C}^\times$ that defines the braiding via~${ c_{\mathbb{C}_g, \mathbb{C}_h} = c(g,h) \id_{\mathbb{C}_{g+h}}}$.
The hexagon equations translate into explicit conditions on the pair~${(\omega,c)}$, which is called an \emph{abelian $3$-cocycle} (see \cite[Section 8.4]{Book-ENO}).

The map $q\colon A \to \mathbb{C}^\times$ defined by $q(l) = c(l,l)$ is a \emph{quadratic form} on $A$, meaning that $q(-l) = q(l)$ for all $l \in A$ and the symmetric map
$b_q(k,l) := q(k+l)q(k)^{-1}q(l)^{-1}$, $ k,l \in A$,
is a bicharacter. By \cite[Theorem 26.1]{EM2}, the quadratic form $q$ completely determines the abelian cohomology class of $(\omega,c)$. Since the quadratic form encodes all the relevant information, we adopt the notation $\operatorname{Vec}_A^q$ for the braided fusion category associated to the pair $(A,q)$.

The group of equivalence classes of braided autoequivalences of $\operatorname{Vec}_A^q$ is naturally isomorphic to $O(A,q)$, the group of automorphisms of $A$ that stabilize $q$. Since $b_q(g,h) = c(g,h)c(h,g)$, the category $\operatorname{Vec}_A^q$ is non-degenerate if and only if the $S$-matrix $S(g,h) = b_q(g,h)$ is invertible, which occurs precisely when $b_q$ is a non-degenerate symmetric bicharacter.
\end{Example}
\subsection{Drinfeld center of a fusion category}

An important class of non-degenerate fusion categories can be constructed using the Drinfeld center $\mathcal{Z}(\cC)$ of a fusion category $(\cC,a,\mathbf{1})$. The center construction produces a non-degenerate braided fusion category $\mathcal{Z}(\cC)$ from any fusion category $\cC$.

Objects of $\mathcal{Z}(\cC)$ are pairs $(Z, \sigma_{-,Z})$, where $Z \in \cC$ and $\sigma_{-,Z} \colon - \otimes Z \to Z \otimes -$ is a natural isomorphism such that the diagram
\begin{equation*}
\xymatrix{
& X \otimes (Z \otimes Y) \ar[r]^{a_{X,Z,Y}^{-1}} & (X \otimes Z) \otimes Y &\\
X \otimes (Y \otimes Z) \ar[ru]^{\id_X \otimes \sigma_{Y,Z}} \ar[rd]_{a_{X,Y,Z}^{-1}} &&& (Z \otimes X) \otimes Y \ar[lu]_{\sigma_{X,Z} \otimes \id_Y}\\
& (X \otimes Y) \otimes Z \ar[r]_{\sigma_{X \otimes Y,Z}} & Z \otimes (X \otimes Y) \ar[ru]_{a_{Z,X,Y}^{-1}} &
}
\end{equation*}
commutes for all $X,Y \in \cC$.

The braided tensor structure on $\mathcal{Z}(\cC)$ is defined as follows:
\begin{itemize}\itemsep=0pt
\item The tensor product is $(Y, \sigma_{-,Y}) \otimes (Z, \sigma_{-,Z}) = (Y\otimes Z, \sigma_{-,Y\otimes Z})$, where the natural isomorphism $\sigma_{X,Y \otimes Z} \colon X \otimes (Y \otimes Z) \to (Y \otimes Z) \otimes X$ is defined by the commutative diagram:
\begin{equation*}
\xymatrix{
X \otimes (Y \otimes Z) \ar[r]^{a_{X,Y,Z}^{-1}} \ar[d]_{\sigma_{X,Y \otimes Z}} & (X \otimes Y) \otimes Z \ar[r]^{\sigma_{X,Y} \otimes \id_Z} & (Y \otimes X) \otimes Z \ar[d]^{a_{Y,X,Z}} \\
(Y \otimes Z) \otimes X \ar[r]_{a_{Y,Z,X}} & Y \otimes (Z \otimes X) \ar[r]_{\id_Y \otimes \sigma_{X,Z}} & Y \otimes (X \otimes Z).
}
\end{equation*}
\item The braiding is given by the isomorphism $\sigma_{X,Y}$.
\end{itemize}

By \cite[Proposition 9.3.4]{Book-ENO}, for any fusion category $\cC$ we have
$\FPdim(\mathcal{Z}(\cC))=\FPdim(\cC)^2$.

For a braided fusion category $\cb$, there is a canonical braided embedding functor
$ \cB \to \mathcal{Z}(\cB)$, $ X\mapsto (X,c_{-,X})$.
This embedding will be important in our analysis of fermionic fusion categories and their minimal non-degenerate extensions.

\subsection{Group actions and (de)equivariantization} \label{categorical actions}

In this subsection, we recall the constructions of equivariantization and de-equivariantization for fusion categories, following \cite{DGNO}.

Let $\cC$ be a fusion category. We denote by $\underline{\operatorname{Aut}_\otimes(\cC)}$ the monoidal category whose objects are tensor autoequivalences of $\cC$, arrows are tensor natural isomorphisms, and tensor product is composition of functors. An \emph{action} of a finite group $G$ on $\cC$ is a monoidal functor $*\colon\underline{G}\to \underline{\operatorname{Aut}_\otimes(\cC)}$, where $\underline{G}$ denotes the discrete monoidal category with objects as elements of $G$.

Such an action provides tensor functors $(g_*,\psi(g))\colon \cC\to \cC$ for each $g \in G$ and natural tensor isomorphisms $\phi(g,h)\colon (gh)_*\to g_*\circ h_*$ for all $g,h \in G$, satisfying coherence conditions \cite[Section 2]{Tam-act}. We call such actions \emph{categorical actions} to distinguish them from the fermionic actions introduced later. We will assume without loss of generality that $(e_*,\psi(e))=\Id_\cC$ is the identity tensor functor and $\phi(e,g)=\phi(g,e)=\id$ are the identity natural transformations \cite[Proposition 3.1]{galindo2017coherence}.

Given fusion categories $\cC$ and $\cD$ with categorical actions of a finite group $G$, a \emph{$G$-equivariant tensor functor} is a pair $(F, \eta)$, where $F \colon \cC \to \cD$ is a tensor functor and $\eta(g) \colon g_* \circ F \to F \circ g_*$ is a family of tensor natural isomorphisms indexed by $G$, such that $\eta(e) = \operatorname{Id}_F$ and for all $X \in \operatorname{Obj}(\cC)$, $g, h \in G$ the diagram
\[
\begin{tikzcd}
(gh)_* F(X) \ar{rr}{\eta(gh)_X} \ar{dd}{\phi(g,h)_{F(X)}}&& F((gh)_*(X)) \ar{dd}{F(\phi(g,h)_{X})} \\\\
g_*h_*(F(X)) \ar{rd}{g_*(\eta(h)_{F(X)})}&& F(g_*h_*(X))\\
&g_*(F(h_*(X))) \ar{ur}{\eta(g)_{h_*(X)}}&
\end{tikzcd}
\]
commutes.

We say that $(F, \eta)$ is an \emph{equivalence of $G$-categories} if the functor $F$ is an equivalence of categories.

If $(F, \eta), (L, \chi) \colon \cC \to \cD$ are $G$-equivariant tensor functors, a \emph{$G$-equivariant tensor natural transformation} $\varphi \colon F \to L$ is a tensor natural transformation such that the diagrams
\[
\begin{tikzcd}
F(g_*(X)) \ar{r}{\varphi_{g_*(X)} } & L(g_*(X)) \\
g_*(F(X)) \ar{u}{\eta(g)_X} \ar{r}{g_*(\varphi_{X})}& g_*(L(X))\ar{u}{\chi(g)_X}
\end{tikzcd}
\]
commute for all $X \in \cC$ and $g \in G$.

The above definitions give rise, for a fixed finite group $G$, to the 2-category of fusion categories with categorical $G$-actions, where objects are fusion categories equipped with a categorical $G$-action, 1-morphisms are $G$-equivariant tensor functors, and 2-morphisms are $G$-equivariant tensor natural transformations.

\subsubsection{The equivariantization construction}

Given a categorical action $*\colon \underline{G}\to \underline{\operatorname{Aut}_\otimes(\cC)}$, the \emph{$G$-equivariantization} $\cC^G$ is the fusion category whose objects are pairs $(V, \tau)$, where $V \in \cC$ and $\tau = \{\tau_g\mid g_*(V) \to V\}_{g \in G}$ satisfies the coherence condition
$
\tau_{gh}= \tau_g \circ g_*(\tau_h) \circ \phi(g,h)_V
$
for all $g,h \in G$. A morphism $\sigma\colon (V, \tau) \to (W, \tau')$ between $G$-equivariant objects is a morphism $\sigma\colon V \to W$ in $\cC$ such that $\tau'_g \circ g_*(\sigma) = \sigma \circ \tau_g$ for all $g \in G$. The tensor structure is given by
\[
(V, \tau)\otimes (W, \tau') = (V\otimes W, \tau''),
\] where
\smash{$
\tau''_g = (\tau_g \otimes \tau'_g) \circ \psi(g)_{V,W}^{-1}$}.

Every $G$-equivariant tensor functor $(F, \eta)\colon \cC \to \cD$ induces a tensor functor $F^G\colon \cC^G \to \cD^G$ between the respective $G$-equivariantizations. If $(V, \tau)$ is an object in $\cC^G$, its image under $F^G$ is the pair
$
F^G(V, \tau) := (F(V), \rho)$,
where $\rho = \{\rho_g\mid g_*(F(V)) \to F(V)\}_{g \in G}$ is defined by the composition
$\rho_g := F(\tau_g) \circ \eta(g)_V$.

Recall that if $\cB$ is a braided fusion category and $\cC$ is a fusion category, a \emph{central functor} from~$\cB$ to $\cC$ is a braided functor $\cB \to \mathcal{Z}(\cC)$. If~$\cE$ is a symmetric fusion category, a \emph{fusion category over} $\cE$ is a fusion category $\cC$ endowed with a \emph{central inclusion} $\cE \to \mathcal{Z}(\cC)$ such that its composition with the forgetful functor is an inclusion $\cE \hookrightarrow \cC$. If $\cC$ is braided, it is a \emph{braided fusion category over} $\cE$ if it is endowed with a braided inclusion $\cE \to \mathcal{Z}_2(\cC)$ (see \cite{DGNO,ENO2}).

The equivariantization $\cC^G$ is canonically a category over $\Rep(G)$. By definition, this means~$\cC^G$ is equipped with a central functor $\Rep(G) \to \mathcal{Z}\bigl(\cC^G\bigr)$ such that composition with the forgetful~functor yields a faithful embedding $\Rep(G) \hookrightarrow \cC^G$. This structure arises from the equiv\-ari\-anti\-zation of the canonical $G$-equivariant tensor functor $I\colon\operatorname{Vec} \to \cC$, where $G$ acts trivially on~$\operatorname{Vec}$, obtaining $I^G\colon \Rep(G)\to \cC^G$.

When the $G$-action consists of braided autoequivalences, $\cC^G$ inherits a braided structure and~$I^G$ becomes a braided functor \cite[Theorem 4.18\,(ii)]{DGNO}.

\subsubsection{De-equivariantization}

Let $\cD$ be a fusion category over $\Rep(G)$, equipped with a central functor $\Rep(G) \to \mathcal{Z}(\cD)$. The \emph{de-equivariantization} $\cD_G$ is defined as the category of modules over the regular algebra $A = \operatorname{Fun}(G)$ acting on $\cD$ via the central functor \cite[Sections 4.1.4 and~4.1.10]{DGNO}.

The category $\cD_G$ is naturally equipped with a canonical $G$-action arising from the right translation action of $G$ on the regular algebra $\operatorname{Fun}(G)$ \cite[Section 4.2.4]{DGNO}. By \cite[Theorem 4.18\,(i)]{DGNO} and \cite[Proposition 4.19]{DGNO}, equivariantization and de-equivariantization are inverse processes that define a 2-equivalence between the 2-category of fusion categories with $G$-actions and fusion categories over $\Rep(G)$.

These constructions relate the Frobenius--Perron dimensions as follows:
\[
 \FPdim\bigl(\cC^G\bigr) = |G| \cdot \FPdim(\cC),\qquad \FPdim(\cD_G) = \FPdim(\cD)/|G|,
\]
when the functor $\Rep(G) \to \cD$ is faithful \cite[Proposition~4.26]{DGNO}.

\subsubsection[Braided G-crossed fusion categories]{Braided $\boldsymbol{G}$-crossed fusion categories}

When working with braided fusion categories, de-equivariantization naturally produces $G$-cros\-sed braided structures. We first recall that a fusion category $\mathcal{C}$ is graded over a finite group $G$, or \emph{$G$-graded}, if it decomposes as a direct sum $\mathcal{C} = \bigoplus_{g \in G} \mathcal{C}_g$ where each component $\mathcal{C}_g$ is a full abelian subcategory and the tensor product respects the grading: $\mathcal{C}_g \otimes \mathcal{C}_h \subseteq \mathcal{C}_{gh}$ for all $g,h\in G$. Objects in $\mathcal{C}_g$ are said to have \emph{degree} $g$. The $G$-grading is called \emph{faithful} if $\mathcal{C}_g \neq 0$ for all $g \in G$.

\begin{Definition}[$G$-crossed braided fusion category]
A \emph{$G$-crossed braided fusion category} $\cB$ is a fusion category equipped with three compatible structures:
\begin{enumerate}\itemsep=0pt
\item[(1)] A $G$-grading $\cB = \bigoplus_{g \in G} \cB_g$.
\item[(2)] An action of the group $G$ on $\cB$ such that $h_*(\cB_g) \subseteq \cB_{hgh^{-1}}$ for all $g,h \in G$.
\item[(3)] Natural isomorphisms $c_{X,Y}\colon X \otimes Y \to g_*(Y) \otimes X$ for $X \in \cB_g$ and $Y \in \cB$, where $g_*(Y)$ denotes the action of $g$ on $Y$ (called the $G$-braiding).
\end{enumerate}
These structures must satisfy compatibility conditions ensuring that the $G$-braiding is functorial in both variables, respects the $G$-action, and satisfies hexagon-type axioms relating the $G$-braiding with associativity. We refer the reader to \cite[Definition 4.41]{DGNO} for the detailed compatibility conditions.
\end{Definition}

The fundamental relationship is that if $\mathcal{D}$ is a braided category over $\Rep(G)$, then its de-equivariantization $\mathcal{D}_G$ naturally inherits the structure of a braided $G$-crossed category \cite[Proposition 4.55\,(i)]{DGNO}. The grading is faithful if and only if the central functor $\Rep(G) \to \mathcal{D}$ is faithful. Conversely, the equivariantization of a braided $G$-crossed category yields a braided category containing $\Rep(G)$. By \cite[Theorem 4.44]{DGNO}, this establishes a 2-equivalence between braided $G$-crossed fusion categories and braided fusion categories containing $\Rep(G)$.

\subsubsection{Centralizers and non-degeneracy}\label{subsection: centralizer}

For a braided fusion category $\mathcal{D}$ over $\Rep(G)$ with de-equivariantization $\mathcal{C} = \mathcal{D}_G$, the centralizer~$C_\mathcal{D}(\Rep(G))$ relates to the trivial component of the $G$-crossed structure on $\mathcal{C}$. Specifically, centralizers commute with de-equivariantization: $(C_\mathcal{D}(\Rep(G)))_G \simeq C_\mathcal{C}(\Rep(G)_G)$ \cite[Proposition~4.30\,(iii)]{DGNO}.

The trivial component $\mathcal{C}_1$ of a braided $G$-crossed category $\mathcal{C} = \mathcal{D}_G$ is given by \[
\mathcal{C}_1 = (C_\mathcal{D}(\Rep(G)))_G,
\]
 and its equivariantization recovers the centralizer: $C_\mathcal{D}(\Rep(G)) = \mathcal{C}_1^G$ \cite[Proposition 4.56\,(i)]{DGNO}.

For non-degeneracy, recall that a braided fusion category $\cB$ is \emph{non-degenerate} if its M\"uger center $\mathcal{Z}_2(\cB)$ is trivial \cite[Theorem 3.4]{DGNO}. The key relationship is that a braided category $\mathcal{D}$ over~$\Rep(G)$ is non-degenerate if and only if the trivial component $\mathcal{C}_1$ of its de-equivarianti\-za\-tion~${\mathcal{C} = \mathcal{D}_G}$ is non-degenerate and the $G$-grading on $\mathcal{C}$ is faithful \cite[Proposition 4.56\,(ii)]{DGNO}. This criterion is used for constructing minimal non-degenerate extensions in Section \ref{section: minimal non-degenerated extension}.

\subsection{Obstruction theory to categorical actions}\label{subsection: obstruction}

To conclude this preliminary section, we recall the construction of the $H^3$-obstruction associated with a group homomorphism $\rho\colon G\to \operatorname{Aut}_\otimes(\cC)$, where $\operatorname{Aut}_\otimes(\cC)$ denotes the group of isomorphism classes of tensor auto-equivalences of $\cC$.

\begin{Definition}
Let $\rho\colon G\to \Autr{\cC}$ be a group homomorphism, where $\cC$ is a fusion category and $G$ is a group. A \emph{lifting} of $\rho$ is a monoidal functor $\tilde{\rho}\colon\underline{G}\to \underline{\Autr{\cC}}$ such that the isomorphism class of $\tilde{\rho}(g)$ is $\rho(g)$ for each $g\in G$.
\end{Definition}

Let $\cC$ be a fusion category and define
\begin{align*}
\widehat{K_0(\mathcal{C})} := \{f\colon \text{Irr}(\mathcal{C})  \to \mathbb{C}^\times \mid f(Y) = f(X_1)f(X_2) \text{ whenever } Y \text{ is a subobject of } X_1 \otimes X_2\}.
\end{align*}

Thus \smash{$\widehat{K_0(\cC)}$} forms an abelian group under pointwise multiplication. Moreover, for every tensor autoequivalence $F \in \operatorname{Aut}_\otimes(\mathcal{C})$, the abelian group $\operatorname{Aut}_\otimes(F)$ of natural tensor automorphisms of $F$ can be canonically identified with \smash{$\widehat{K_0(\mathcal{C})}$} as follows: if $\Phi \in \operatorname{Aut}_\otimes(F)$, then for every simple object $X \in \text{Irr}(\mathcal{C})$, the component $\Phi_X\colon F(X) \to F(X)$ is given by $\Phi_X = f(X) \cdot \text{id}_{F(X)}$ for some scalar $f(X) \in \mathbb{C}^\times$, and this assignment defines a function \smash{$f \in \widehat{K_0(\mathcal{C})}$}.

Let $\rho\colon G \to \operatorname{Aut}_\otimes(\mathcal{C})$ be a group homomorphism. Then $G$ acts on $\operatorname{Irr}(\mathcal{C})$ in a way that induces algebra automorphisms on $K_0(\mathcal{C})$. This induces a natural $G$-module structure on~\smash{$\widehat{K_0(\mathcal{C})}$} given~by
\[
(g \cdot f)(X) = f\bigl(\rho\bigl(g^{-1}\bigr)(X)\bigr)
\] for $g \in G$, \smash{$f \in \widehat{K_0(\mathcal{C})}$}, and $X \in \operatorname{Irr}(\mathcal{C})$. We denote by~\smash{$H^n_\rho\bigl(G, \widehat{K_0(\mathcal{C})}\bigr)$} the $n$-th cohomology group of $G$ with coefficients in this $G$-module structure induced by $\rho$.

Let us fix a representative tensor autoequivalence $g_*\colon \mathcal{C} \to \mathcal{C}$ for each $g \in G$ and a tensor natural isomorphism $\theta_{g,h}\colon g_* \circ h_* \to (gh)_*$ for each pair $g,h \in G$. The failure of these choices to define a categorical action is measured by the obstruction \smash{$O_3(\rho)(g,h,l) \in \widehat{K_0(\mathcal{C})}$} defined by the commutativity of the diagram
\begin{equation}\label{diagram: obstruction 3}
\begin{gathered}
\xymatrixcolsep{5pc} \xymatrix{
g_*\circ h_*\circ l_* \ar[dd]^{g_*(\theta_{h,l})} \ar[r]^{(\theta_{g,h})_{l_*}}& (gh)_*\circ l_* \ar[d]^{\theta_{gh,l}}\\
 &(ghl)_*\ar@{..>}[d]^{O_3(\rho)(g,h,l)}\\
g_*\circ (hl)_* \ar[r]^{\theta_{g,hl}} &(ghl)_*.
 }
\end{gathered}
\end{equation}

\begin{Proposition}[{\cite[Theorem 5.5]{Ga1}}]\label{prop:3-cocycle-obstruction}
Let $\cC$ be a fusion category and $\rho\colon G\to \operatorname{Aut}_\otimes(\cC)$ a group homomorphism. The map \smash{$O_3(\rho,\theta)\colon G^{\times 3}\to \widehat{K_0(\C)}$} defined by diagram \eqref{diagram: obstruction 3} is a $3$-cocycle, and its cohomology class \smash{$O_3(\rho)\in H_{\rho}^3\bigl(G,\widehat{K_0(\C)}\bigr) $} depends only on $\rho$. The homomorphism $\rho$ lifts to an action \smash{$\widetilde{\rho}\colon\underline{G}\to \underline{\operatorname{Aut}_\otimes(\cC)}$} if and only if $O_3(\rho)=0$.
\end{Proposition}

When $O_3(\rho)$ vanishes, categorical liftings exist and their non-uniqueness is controlled by \smash{$H_{\rho}^2\bigl(G,\widehat{K_0(\C)}\bigr)$}. Specifically, given a lifting \smash{$*\colon\underline{G}\to \underline{\operatorname{Aut}_\otimes(\cC)}$} and a 2-cocycle \smash{$\beta \in Z^2_\rho\bigl(G, \widehat{K_0(\cC)}\bigr)$}, we can construct a new lifting via the commutativity of the diagrams

\begin{equation}\label{eq: torciendo accion}
\xymatrix{
(gh)_* \ar@{..>}[rr]^{(\beta \triangleright \phi)(g,h)} \ar[dr]_{\beta(g,h)} && g_*\circ h_* \\
& g_* \circ h_*\ar[ur]_{\phi(g,h)}. &
}
\end{equation}

The data $(g_*, \psi(g), (\beta \triangleright \phi)(g,h))$ defines a new lifting. More precisely, we have the following result:

\begin{Proposition}[{\cite[Theorem 5.5]{Ga1}}]\label{Proposition:Obstruction-bosonic}
Let $\cC$ be a fusion category and $\rho\colon G\to \operatorname{Aut}_\otimes(\cC)$ a~group homomorphism. If $[O_3(\rho)]=0$, the set of equivalence classes of liftings of $\rho$ is a torsor over~\smash{$H^2_\rho\bigl(G,\widehat{K_0(\C)}\bigr)$}, where the torsor action is given by the twisting construction $\beta \triangleright \phi$ described in \eqref{eq: torciendo accion}.
\end{Proposition}

\section{Actions of super-groups on fermionic fusion categories}\label{acciosuper}

In this section, we develop the theory of fermionic actions on fermionic fusion categories and establish their relationship with fusion categories over super-groups. We first introduce fermionic fusion categories and their morphisms, then define and study fermionic actions of super-groups.

To avoid lengthy definitions involving unit isomorphisms while maintaining generality, we adopt the following conventions throughout this paper.

All monoidal categories are assumed to be strict with respect to the unit object, meaning~${V \otimes \mathbf{1} = V = \mathbf{1} \otimes V}$ and $a_{V,\mathbf{1},W} = a_{\mathbf{1},V,W} = \id_{V \otimes W}$ for all objects $V$, $W$. Similarly, for monoidal functors $(F,\tau)\colon \cC \to \cD$, we assume $F(\mathbf{1}) = \mathbf{1}$ and $\tau_{\mathbf{1},V} = \tau_{V,\mathbf{1}} = \id_{F(V)}$ for all $V \in \cC$.

By \cite[Theorem 3.2]{schauenburg2001turning}, any monoidal category is monoidally equivalent to one with strict unit without changing the underlying category, only slightly redefining the tensor product. In particular, if the original category is skeletal, the unit strictified category remains skeletal. Similarly, by \cite[Proposition 3.1]{galindo2017coherence}, any strong monoidal functor can be made strict with respect to the unit. These results justify our assumptions without loss of generality.

\subsection{The 2-category of fermionic fusion categories}

We begin by defining the 2-category of fer\-mionic fusion categories, which is by definition the 2-category of fusion categories over the symmetric category of super vector spaces.

\begin{Definition}\label{definition: spin y fermionic categories}
Let $\cC$ be a fusion category. An object $(f,\sigma_{-,f})\in \mathcal{Z}(\cC)$ is called a \emph{fermion} if~${f \otimes f \cong \mathbf{1}}$ and $\sigma_{f,f}= -\id_{f \otimes f}$.
\begin{enumerate}[leftmargin=*,label={\rm (\alph*)}]\itemsep=0pt
 \item A \emph{fermionic fusion category} is a pair $(\cC,(f,\sigma_{-,f}))$ consisting of a fusion category $\cC$ and a~fermion $(f,\sigma_{-,f}) \in \mathcal{Z}(\cC)$.
 \item A \emph{spin-braided fusion category} is a pair $(\cB,f)$ where $\cB$ is a braided fusion category and $(f,c_{-,f}) \in \mathcal{Z}(\cB)$ is a fermion.
\end{enumerate}
\end{Definition}

\begin{Remark}\label{remark:fermion-properties}\qquad
\begin{enumerate}[leftmargin=*,label={\rm (\alph*)}]\itemsep=0pt
\item If $(f,\sigma_{-,f})$ is a fermion, then $f$ is not isomorphic to the unit object $\mathbf{1}$. This follows from the definition of the half-braiding, as $\sigma_{\mathbf{1},\mathbf{1}}= \id_\mathbf{1}$.

\item If $(f,\sigma_{-,f})\in \mathcal{Z}(\cC)$ is a fermion, the fusion subcategory $\langle (f,\sigma_{-,f}) \rangle$ of $\mathcal{Z}(\cC)$ is braided equivalent to $\SV$. Therefore, a fermionic fusion category $(\cC,(f,\sigma_{-,f}))$ is precisely a fusion category over $\SV$ in the sense of \cite[Definition 4.16]{DGNO}.
\item Fermionic structures on a fusion category $\mathcal{C}$ correspond precisely to spin-braided structures on its Drinfeld center $\mathcal{Z}(\mathcal{C})$. Thus, classifications of spin-braided structures on centers are equivalent to classifications of fermionic structures on the underlying categories.
\end{enumerate}
\end{Remark}

We present several examples of fermionic fusion categories that will be important throughout this paper.

\begin{Example}[fermionic pointed fusion categories] \label{fermion pointed fusion categories}
According to \cite[Proposition 2.6]{Fermioni-modular}, there is a correspondence between fermions in $\operatorname{Vec}_G^\omega$ and pairs $(f,\eta)$, where $\eta\colon G\to \mathbb{C}^\times$ satisfies:
\begin{enumerate}[leftmargin=*,label={\rm (\alph*)}]\itemsep=0pt
 \item $f\in Z(G)$ is of order two,
 \item $\frac{\eta(x)\eta(y)}{\eta(xy)}=\frac{\omega(f,x,y)\omega(x,y,f)}{\omega(x,f,y)}$ for all $x, y \in G$,
 \item $\frac{\omega(x,f,f)\omega(f,x,f)}{\omega(f,f,x)}\eta(x)^2=1$ for all $x\in G$,
 \item $\eta(f)=-1$.
\end{enumerate}
This correspondence follows from the equivalence between fermionic structures on $\operatorname{Vec}_G^\omega$ and spin-braided structures on $\mathcal{Z}\bigl(\operatorname{Vec}_G^\omega\bigr)$ mentioned in Remark \ref{remark:fermion-properties}\,(c).
\end{Example}

\begin{Example}[Ising categories as spin-braided fusion categories]\label{Example:Ising}
An Ising fusion category is a non-pointed fusion category of Frobenius--Perron dimension 4 (see \cite[Appendix B]{DGNO}). These categories have three simple objects $\mathbf{1}$, $f$, $\sigma$, with fusion rules
$
\sigma^2=\mathbf{1}+f$, $ f^2=\mathbf{1}$, $ f\sigma=\sigma f=\sigma$.

The associativity constraints are given by the $F$-matrices
\[
 F_{\sigma\sigma\sigma}^\sigma=\frac{\epsilon}{\sqrt{2}}\left(\begin{matrix}
1 & 1 \\
1 & -1 \end{matrix}\right), \qquad F_{f \sigma f}^\sigma=F_{\sigma f\sigma}^f=-1,
\]
where $\epsilon\in \{1,-1\}$.

Ising fusion categories admit braided structures, and in all cases $f$ is a fermion. An Ising braided fusion category is always non-degenerate and $\mathcal{I}_{\rm pt}= \langle f \rangle \cong \SV$. As proved in \cite[Appendix B]{DGNO}, there are 8 equivalence classes of Ising braided fusion categories.
\end{Example}

To simplify definitions and avoid working with a fixed isomorphism $u\colon f \otimes f \to \mathbf{1}$, we will assume from now on that $f \otimes f = \mathbf{1}$. This represents no loss of generality since, by the same technique as in \cite[Theorem 3.2]{schauenburg2001turning}, we can slightly modify the tensor product without changing the underlying category to obtain an equivalent one with the property $f \otimes f = \mathbf{1}$. In particular, if we are working with a skeletal category, the modified category remains skeletal.

\begin{Definition}\label{fermionicfunctor}
Let $(\cC,(f,\sigma_{-,f}))$ and $(\cC',(f',\sigma'_{-,f'}))$ be fermionic fusion categories. A \emph{fermionic tensor functor} is a pair $((F,\tau), \phi)$, where $(F,\tau)\colon \cC\to \cC'$ is a tensor functor and $\phi\colon F(f)\to f'$ is an isomorphism satisfying the following conditions:
\begin{enumerate}\itemsep=0pt
\item[(1)]
\begin{equation}\label{eq: fermionic funtor condicion 1}
 \tau_{f,f}^{-1} = \phi \otimes \phi\colon \ F(f)\otimes F(f)\to \unit.
\end{equation}
\item[(2)] The compatibility diagram commutes for each $V\in \cC$
\begin{equation}\label{diagram:fermion functor}
\xymatrix{
F(V\otimes f) \ar[d]_{\tau_{V,f}} \ar[rr]^{F(\sigma_{V,f})} && F(f\otimes V)\ar[d]^{\tau_{f,V}} \\
F(V)\otimes F(f)\ar[d]_{\id_{F(V)} \otimes \phi} && F(f)\otimes F(V) \ar[d]^{\phi \otimes \id_{F(V)}} \\
F(V) \otimes f' \ar[rr]_{\sigma'_{F(V),f'}} && f' \otimes F(V).
}
\end{equation}
\end{enumerate}
\end{Definition}
The following proposition records basic properties of fermionic tensor functors that will be used throughout this paper.

\begin{Proposition}\label{prop:fermionic-functor-properties}
Let $F\colon \cC \to \cC'$ be a tensor functor between fermionic fusion categories such that $F(f) \cong f'$. Then:
\begin{enumerate}\itemsep=0pt
\item[$(1)$] There always exists an isomorphism $\phi\colon F(f) \to f'$ satisfying condition \eqref{eq: fermionic funtor condicion 1}.
\item[$(2)$] The commutativity of diagram~\eqref{diagram:fermion functor} is independent of the choice of $\phi$.
\item[$(3)$] Each tensor functor $F$ admits at most two fermionic structures, corresponding to $\pm\phi$.
\end{enumerate}
\end{Proposition}

\begin{proof}
(1) Take any $\omega\colon F(f) \to f'$. This defines $c \in \mathbb{C}^\times$ by $c \cdot \id_{\unit'} = (\omega \otimes \omega) \circ \tau_{f,f}$, then $\phi = \frac{1}{\sqrt{c}} \cdot \omega$.

(2) and (3) follows immediately from the proof of (1).
\end{proof}

\begin{Example}
If $(\cB,c)$ and $(\cB',c')$ are spin-braided fusion categories with fermions $f$ and $f'$ respectively, and $F\colon\cB\to \cB'$ is a braided functor such that $\phi\colon F(f)\cong f'$, then $F$ satisfies diagram \eqref{diagram:fermion functor} for any isomorphism $\phi\colon F(f)\to f'$.
\end{Example}

\begin{Example}[Fermionic functors between pointed categories]\label{Example:Fermionic functors in the pointed case}
Let \smash{$\bigl(\operatorname{Vec}_{G_1}^{\omega_1},(f,\eta)\bigr)$} and $\smash{\bigl(\operatorname{Vec}_{G_2}^{\omega_2}},\allowbreak(f',\eta')\bigr)$ be fermionic pointed fusion categories as in Example \ref{fermion pointed fusion categories}. A fermionic tensor functor between these categories consists of the following data:
\begin{enumerate}[leftmargin=*,label={\rm (\alph*)}]\itemsep=0pt
 \item A group homomorphism $F\colon G_1 \to G_2$ with $F(f) = f'$.
 \item A scalar $c \in \mathbb{C}^\times$.
 \item A normalized 2-cochain $\tau\colon G_1^{\times 2} \to \mathbb{C}^\times$.
\end{enumerate}
These data must satisfy the following compatibility conditions:
\begin{enumerate}[leftmargin=*,label={\rm (\roman*)}]\itemsep=0pt
\item $
 \omega_1(g,h,l)\tau(gh,l)\tau(g,h) = \tau(g,hl)\tau(h,l)\omega_2(F(g),F(h),F(l))
$
 for all $g,h,l \in G_1$.
 \item $\tau(f,f) = c^2$.
 \item
$
 \frac{\tau(f,g)}{\tau(g,f)} = \frac{\eta'(F(g))}{\eta(g)}
$
 for all $g \in G_1$.
\end{enumerate}
The functor is defined by \smash{$F_*(\mathbb{C}_g) = \mathbb{C}_{F(g)}$}, \smash{$\tau_{\mathbb{C}_g,\mathbb{C}_h} = \tau(g,h)\id_{\mathbb{C}_{F(gh)}}$}, and \smash{$\phi = c^{-1}\id_{\mathbb{C}_{f'}}$}.
\end{Example}

\begin{Definition}
Let $((F,\tau),\phi), ((F',\tau'),\phi')\colon\cC\to \cD$ be fermionic tensor functors. A tensor natural transformation $\gamma\colon F\to F'$ is called a \emph{fermionic tensor natural transformation} if the diagram
\begin{equation}\label{eq: diagram fermionic natural}
\begin{tikzcd}
F(f_{\cC}) \arrow[rr, "\gamma_{f_{\cC}}"] \arrow[dr, "\phi"'] && F'(f_{\cC}) \arrow[dl, "\phi'"] \\
& f_{\cD}&
\end{tikzcd}
\end{equation}
commutes.
\end{Definition}

To complete the definition of the 2-category of fermionic fusion categories, we describe composition of 1-morphisms and 2-morphisms. If $(F,\tau,\phi)\colon \cC \to \cD$ and $(G,\rho,\omega)\colon \cD \to \mathcal{E}$ are fermionic tensor functors, their composite is fermionic with isomorphism
$
\omega \circ G(\phi)\colon G(F(f_\cC)) \to f_{\mathcal{E}}$.
The vertical and horizontal composites of fermionic tensor natural transformations are again fermionic natural transformations. This defines the 2-category of fermionic fusion categories.

\subsection[The 2-category of fermionic fusion categories with a G-action]{The 2-category of fermionic fusion categories with a $\boldsymbol{G}$-action}

Let $(\cC,(f, \sigma_{-,f}))$ be a fermionic fusion category. We denote by \smash{$\underline{\operatorname{Aut}_\otimes^{\rm ferm}(\cC,f)}$} the monoidal category whose objects are \emph{fermionic tensor autoequivalences} of $(\cC,(f, \sigma_{-,f}))$ and morphisms are \emph{fermionic tensor natural isomorphisms}. Note that we suppress $\sigma_{-,f}$ from the notation, but it is understood that the half-braiding is part of the fermionic structure.

Given a group $G$, we can define the \emph{$2$-category of fermionic $G$-actions} analogously to how one defines the 2-category of categorical $G$-actions on fusion categories, but replacing fusion categories with fermionic fusion categories throughout. More precisely, this 2-category consists of the following data:
\begin{itemize}\itemsep=0pt
\item \textit{Objects:} Fermionic fusion categories $(\cC, f)$ equipped with a \emph{fermionic $G$-action}, that is, a monoidal functor \smash{$*\colon \underline{G}\to \underline{\operatorname{Aut}_\otimes^{\rm ferm}(\cC,f)}$}.
\item \textit{$1$-morphisms:} \emph{$G$-equivariant fermionic tensor functors}, that is, a $G$-equivariant tensor functor $(F,\eta)$, where both the tensor functor $F$ and the natural isomorphism $\eta$ are fermionic.
\item \textit{$2$-morphisms:} \emph{$G$-equivariant fermionic tensor natural transformations}, that is, tensor natural transformations that are both $G$-equivariant and fermionic.
\end{itemize}

For the precise definitions of the 2-category of fusion categories with categorical $G$-actions, see Section~\ref{categorical actions}.

There is a canonical \emph{forgetful functor}
\smash{$
U\colon\underline{\operatorname{Aut}_\otimes^{\rm ferm}(\cC,f)}\to \underline{\operatorname{Aut}_\otimes(\cC)}
$}
which forgets the fermionic structure. This functor is faithful but not full, since for any fermionic tensor natural automorphism $\gamma\colon F\to F$ we always have $\gamma_f=\operatorname{id}_f$ by the commutative diagram~\eqref{eq: diagram fermionic natural}.

We denote by $\underline{\operatorname{Aut}_\otimes(\cC,f)}\subset \underline{\operatorname{Aut}_\otimes(\cC)}$ the \emph{image} of the forgetful functor $U$. This subcategory consists of those tensor autoequivalences $F\colon \cC\to \cC$ such that $F(f)\cong f$ and diagram \eqref{diagram:fermion functor} commutes for some (equivalently, any) choice of isomorphism $\phi\colon F(f)\to f$. As with previous notation, this subcategory depends on the half-braiding $\sigma_{-,f}$, but we suppress it from the notation for simplicity.

By Proposition \ref{prop:fermionic-functor-properties}, the condition for $F$ to belong to $\underline{\operatorname{Aut}_\otimes(\cC,f)}$ depends only on the tensor functor $F$ and the half-braiding $\sigma_{-,f}$, not on the specific choice of isomorphism $\phi$.

In the braided case, \smash{$\underline{\operatorname{Aut}_\otimes^{\rm br}(\cC,f)}$} consists precisely of those braided tensor autoequivalences that preserve the fermion $f$. We denote by $\operatorname{Aut}_\otimes^{\rm br}(\cC,f)$ the group of isomorphism classes of such autoequivalences.

\begin{Example}[pointed spin-braided categories of rank 4]\label{ex: pointed 4 rank}
We identify roots of unity in $\mathbb{C}$ with $\mathbb{Q}/\mathbb{Z}$. A \emph{spin-braided structure} on a pointed braided fusion category $\operatorname{Vec}_A^q$ associated with a~quadratic form $q\colon A\to \mathbb{Q}/\mathbb{Z}$ over an abelian group $A$ is determined by choosing an element~${f \in A}$ of order $2$ such that $q(f) = \frac{1}{2}$. The group $\operatorname{Aut}(A)$ acts naturally on the set of spin-braided structures via $\bigl(F^*q, F^{-1}(f)\bigr)$. The stabilizer of the pair $(q,f)$ under this action is~${O(A,q,f) := \{F \in \operatorname{Aut}(A) \mid F^*q = q \text{ and } F(f) = f\} \subset O(A,q)}$.
The group $O(A,q,f)$ coincides with $\operatorname{Aut}_\otimes^{\rm br}\bigl(\operatorname{Vec}_A^q,f\bigr)$, while $O(A,q)$ is $\operatorname{Aut}_\otimes^{\rm br}\bigl(\operatorname{Vec}_A^q\bigr)$.

We describe the classification of pointed non-degenerate spin-braided fusion categories of rank~4 following \cite[Lemma~A.11]{DGNO}, up to conjugation by $\operatorname{Aut}(A)$.

\textit{Case $1$:} $A = \mathbb{Z}/2\mathbb{Z} \times \mathbb{Z}/2\mathbb{Z}$. The non-degenerate quadratic forms are
\begin{gather*}
q_0((x,y)) = \frac{xy}{2}, \qquad q_2((x,y)) = \frac{x^2+y^2}{4},\\
q_{4}((x,y)) = \frac{x^2+xy+y^2}{2},\qquad
 q_6((x,y)) = -\frac{x^2+y^2}{4}.
\end{gather*}
In all cases, $f = (1,1)$ are fermions and $O(A,q_k,f) = \{\operatorname{id}, L\} \cong \mathbb{Z}/2\mathbb{Z}$, where $L(x,y) = (y,x)$.

\textit{Case $2$:} $A = \mathbb{Z}/4\mathbb{Z}$. The non-degenerate quadratic forms are
\smash{$
q_{k}(\bar{x}) = \frac{kx^2}{8}$},
where $k\in \{1,3,5,7\}$.
In all cases, $f = \bar{2}$ and $O(A,q_{\pm},f) = \{\operatorname{id}, \iota\} \cong \mathbb{Z}/2\mathbb{Z}$ where $\iota(\bar{x}) = -\bar{x}$.
\end{Example}

\subsection{The 2-category of fermionic actions of a super-group}

Our goal in this subsection is to define fermionic actions of super-groups on fermionic fusion categories.

We recall that a finite \emph{super-group} is a pair $\bigl(\widetilde{G},z\bigr)$, where $\widetilde{G}$ is a finite group and ${z\in \widetilde{G}}$ is a~central element of order two. Every super-group \smash{$\bigl(\widetilde{G},z\bigr)$} gives rise to an exact sequence
\[
1 \longrightarrow \langle z\rangle \longrightarrow \widetilde{G} \longrightarrow \widetilde{G}/\langle z \rangle \longrightarrow 1,
\]
which is determined up to equivalence by a unique element $\alpha \in H^2\bigl(\widetilde{G}/\langle z \rangle ,\mathbb{Z}/2\mathbb{Z}\bigr)$. This cohomological classification allows us to identify a super-group \smash{$\bigl(\widetilde{G},z\bigr)$} with the associated pair $(G,\alpha)$, where \smash{$G:=\widetilde{G}/\langle z \rangle$} and $\alpha \in H^2(G,\mathbb{Z}/2\mathbb{Z})$.

Every super-group induces an action on the category of super-vector spaces $\SV$. Given a~representative 2-cocycle $\alpha\colon G\times G\to \mathbb{Z}/2\mathbb{Z}$ for the central extension, we can explicitly describe this action as follows. The action is defined by taking $g_* = \Id_{\SV}$ for all $g\in G$, with natural isomorphisms given by
\begin{equation}\label{eq:sv-action}
\phi(g,h)_f = (-1)^{\alpha(g,h)} \cdot \id_f,
\end{equation}
where $f\in \SV$ is the fermion. The equivalence class of this $G$-action depends only on the cohomology class of $\alpha$, and we have $\SV^G \cong \Rep\bigl(\widetilde{G},z\bigr)$ as symmetric fusion categories.

To define fermionic actions of super-groups, we use the notion of co-slice 2-categories. Recall that for a 2-category $\mathcal{C}$ and an object $S \in \mathcal{C}$, the \emph{co-slice $2$-category} $\mathcal{C}_{S/}$ has
\begin{itemize}\itemsep=0pt
 \item \textit{Objects:} the 1-morphisms $a\colon S \to A$ in $\mathcal{C}$;
 \item \textit{$1$-morphisms:} from $a\colon S \to A$ to $b\colon S \to B$, the pairs $(f, \phi)$ where $f\colon A \to B$ is a~1-morphism in $\mathcal{C}$ and $\phi\colon fa \cong b$ is a 2-isomorphism in $\mathcal{C}$;
 \item \textit{$2$-morphisms:} from $(f, \phi)$ to $(g, \psi)$, the 2-morphisms $\xi\colon f \to g$ such that $ \phi\cdot (\xi \id_a) = \psi$.
\end{itemize}

We now apply this construction to define fermionic super-group actions.

\begin{Definition}\label{def:fermionic-super-group-action}
Let $(G,\alpha)$ be a super-group. We define the \emph{$2$-category of fermionic fusion categories with $(G,\alpha)$-actions} as the co-slice 2-category under $\SV$ in the 2-category of fermionic $G$-categories, where $\SV$ is equipped with the $G$-action determined by $\alpha \in Z^2(G,\mathbb{Z}/2\mathbb{Z})$ as in~\eqref{eq:sv-action}.
\end{Definition}

The structure of this co-slice 2-category translates into the following concrete definitions:

\begin{Definition}\label{def:fermionic-action}
A \emph{fermionic action} of $(G,\alpha)$ on $(\cC,f)$ consists of a fermionic categorical $G$-action $(g_*,\psi(g),\phi(g,h),\eta_g)$ such that the following diagram commutes
\[
 \begin{tikzcd}
(gh)_*(f) \ar{rr}{\eta_{gh}} \ar{dd}[swap]{\phi(g,h)_{f}}&& f \ar{dd}{(-1)^{\alpha(g,h)}\id_f} \\\\
g_*h_*(f) \ar{rd}[swap]{g_*(\eta_h)}&& f\\
&g_*(f) \ar{ur}[swap]{\eta_g}.&
\end{tikzcd}
\]
\end{Definition}

\begin{Example}[fermionic actions on pointed fusion categories]\label{fermionic action on pointed fusion categories}
Let $(G,\alpha)$ be a super-group and $\bigl(\operatorname{Vec}_F^{\omega},(f,\eta)\bigr)$ a fermionic pointed fusion category. A fermionic action of $(G,\alpha)$ on $\bigl(\operatorname{Vec}_F^{\omega},\allowbreak(f,\eta)\bigr)$ is determined by
\begin{itemize}\itemsep=0pt
 \item A group homomorphism $*\colon G \to \operatorname{Aut}(F)$ specifying how $G$ acts on the group $F$, with the condition $g_*(f) = f$ for all $g \in G$.
 \item Normalized maps $\mu\colon G \times F \times F \to \mathbb{C}^\times$, $\gamma\colon G \times G \times F \to \mathbb{C}^\times$, and $\xi\colon G \to \mathbb{C}^\times$,
satisfying the following conditions:
\begin{enumerate}[leftmargin=*,label={\rm (\alph*)}]\itemsep=0pt
 \item The map $\mu$ is compatible with the associativity constraint
 \[\frac{\omega(a,b,c)}{\omega(g_*(a),g_*(b),g_*(c))} = \frac{\mu(g;b,c)\mu(g;a,bc)}{\mu(g;ab,c)\mu(g;a,b)}.\]

 \item The maps $\mu$ and $\gamma$ satisfy the coherence conditions
 \begin{gather*}\frac{\mu(g;h_*(a),h_*(b))\mu(h;a,b)}{\mu(gh;a,b)} = \frac{\gamma(g,h;ab)}{\gamma(g,h;a)\gamma(g,h;b)},\\
 \gamma(gh,k;a)\gamma(g,h;k_*(a)) = \gamma(h,k;a)\gamma(g,hk;a).\end{gather*}

 \item The fermionic compatibility condition
 \begin{align*}
 \gamma(g,h;f) = (-1)^{\alpha(g,h)}\frac{\xi(g)\xi(h)}{\xi(gh)}, \qquad \mu(g;f,f)=\xi(g)^{-2},
 \end{align*}
 for all $g,h \in G$.
\end{enumerate}
\end{itemize}

The action is implemented as follows: for each $g \in G$, the tensor functor $g_*$ maps $\mathbb{C}_a$ to $\mathbb{C}_{g_*(a)}$, with monoidal structure $\psi(g)_{\mathbb{C}_a,\mathbb{C}_b} = \mu(g;a,b)\operatorname{id}_{\mathbb{C}_{ab}}$ and natural isomorphisms $\phi(g,h)_{\mathbb{C}_{a}} = \gamma(g,h;a)\operatorname{id}_{\mathbb{C}_{(gh)_*(a)}}$ for all $g,h \in G$ and $a \in F$. The isomorphism $\eta_g\colon g_*(\mathbb{C}_f) \to \mathbb{C}_f$ is given by $\eta_g = \xi(g)\operatorname{id}_{\mathbb{C}_f}$ for all $g \in G$.

For spin-braided pointed fusion categories \smash{$\cB = \operatorname{Vec}_A^{(\omega,c)}$}, we impose the additional condition that $G$ acts on $\cB$ by braided tensor automorphisms. This translates to the compatibility condition
$\frac{c(a,b)}{c(g_*(a),g_*(b))} = \frac{\mu(g;g_*(a),g_*(b))}{\mu(g;a,b)}$
for all $g \in G$ and $a,b \in A$.
\end{Example}

\begin{Definition}
If $\cC$ and $\cD$ are fermionic fusion categories with fermionic $(G,\alpha)$-actions, a~fermionic \emph{$(G,\alpha)$-equivariant functor} is a fermionic $G$-equivariant tensor functor $(F,\kappa)\colon \cC\to \cD$ and an isomorphism $\tau\colon F(f_\cC)\to f_\cD$ such that the diagram
\[
\begin{tikzcd}
F(g_*(f_\cC)) \ar{d}[swap]{\kappa_f}\ar{rr}{F(\eta(g))} && F(f) \ar{d}{\tau}\\
g_{*'}(F(f)) \ar{rr}{g_{*'}(\tau)}& & g_*'{f}
\end{tikzcd}
\]
commute for all $g \in G$.
\end{Definition}

\begin{Definition}
If $(F,\kappa,\tau), (F',\kappa',\tau')\colon \cC\to \cD$ are fermionic $(G,\alpha)$-equivariant tensor functors, a fermionic $(G,\alpha)$-equivariant tensor natural transformation is fermionic tensor natural transformation $\omega\colon (F,\kappa)\to (F',\kappa')$ such that the diagram
\[\begin{tikzcd}
 F(f) \ar[swap]{dr}{\tau} \ar{rr}{\omega_f} && F'(f) \ar{dl}{\tau'}\\
 &f&
\end{tikzcd}\]commutes.
\end{Definition}

We now establish the relationship between fermionic fusion categories with fermionic actions and fusion categories over super-groups.

\begin{Theorem}\label{theorem: equivalencia dequivariantizacion- equivariantizacion}
Let $\bigl(\widetilde{G},z\bigr)$ be a finite super-group and $G:=\widetilde{G}/\langle z \rangle$. $G$-equivariantization and $G$-de-equivariantization processes define a biequivalence of $2$-categories between:
\begin{enumerate}\itemsep=0pt
 \item[$(1)$] Fermionic fusion categories with a fermionic action of $\bigl(\widetilde{G},z\bigr)$.
 \item[$(2)$] Fusion categories over $\Rep\bigl(\widetilde{G},z\bigr)$.
\end{enumerate}
\end{Theorem}

\begin{proof}

Let $\bigl(\widetilde{G},z\bigr)$ be a finite super-group and $G:=\widetilde{G}/\langle z \rangle$. The projection $\widetilde{G}\to G$ defines an inclusion \smash{$\Rep(G) \subset \Rep\bigl(\widetilde{G},z\bigr)$} of symmetric fusion categories.

By \cite[Theorem 4.18, Proposition 4.19]{DGNO}, equivariantization and de-equi\-varian\-tization are mutually inverse processes and define a biequivalence of 2-categories between fusion categories with a $G$-action and fusion categories over $\Rep(G)$.

Using de-equivariantization with respect to $G$, we obtain that \smash{$\Rep\bigl(\widetilde{G},z\bigr)_G$} is braided equivalent to $\SV$, so $\SV^G$ is equivalent to $\Rep\bigl(\widetilde{G},z\bigr)$. Moreover, if $(\cC,(f,\sigma_{-,f}))$ is a fermionic fusion category, then we have a central functor $\SV \to \mathcal{Z}(\cC)$ that induces a central functor~${\SV^G \cong \Rep\bigl(\widetilde{G},z\bigr) \to \mathcal{Z}\bigl(\cC^G\bigr)}$. Hence, $G$-equivariantization is a 2-functor from the 2-catego\-ry of fermionic fusion categories with a fermionic action of \smash{$\bigl(\widetilde{G},z\bigr)$} to the 2-category of fusion categories over $\Rep\bigl(\widetilde{G},z\bigr)$.

Conversely, if $\cD$ is a fusion category over $\Rep\bigl(\widetilde{G},z\bigr)$, then $\Rep(G)\subset \Rep\bigl(\widetilde{G},z\bigr)\subset \mathcal{Z}(\cD)$. It follows from \cite[Proposition 2.10]{ENO2} that $\cD$ is a $G$-equivariantization of some fusion category~$\cD_G$ and the $G$-de-equivariantization of \smash{$C_{\mathcal{Z}(\cD)}(\Rep(G))$} is braided equivalent to $\mathcal{Z}(\cD_G)$. Taking $G$-de-equivariantization of the sequence of inclusions
\[
\Rep(G) \subset \Rep\bigl(\widetilde{G},z\bigr)\subset C_{\mathcal{Z}(\cD)}(\Rep(G)),
\]
 we obtain $\Vc \subset \SV \subset \mathcal{Z}(\cD_G)$. Hence $\cD_G$ is a fermionic fusion category. This proves that de-equivariantization is a functor from the 2-category of fusion categories over $\Rep\bigl(\widetilde{G},z\bigr)$ to the 2-category of fermionic fusion categories with fermionic action of \smash{$\bigl(\widetilde{G},z\bigr)$}.

These functors are mutually inverse by \cite[Theorem 4.18 and Proposition 4.19]{DGNO}.
\end{proof}

\begin{Corollary}\label{coro: equivalencia dequivariantizacion- equivariantizacion braided case}
Let $\bigl(\widetilde{G},z\bigr)$ be a finite super-group and $G:=\widetilde{G}/\langle z \rangle$. $G$-equivariantization and $G$-de-equivariantization processes define a biequivalence of $2$-categories between:
\begin{enumerate}\itemsep=0pt
 \item[$(1)$] Spin-braided fusion categories with a fermionic action of $\bigl(\widetilde{G},z\bigr)$ compatible with the braiding.
 \item[$(2)$] Braided fusion categories $\cD$ over $\Rep\bigl(\widetilde{G},z\bigr)$ such that $\Rep(G)\subseteq \mathcal{Z}_2(\cD)$.
\end{enumerate}
\end{Corollary}

\begin{proof}
This follows directly from Theorem~\ref{theorem: equivalencia dequivariantizacion- equivariantizacion} and \cite[Theorem 4.18 and Proposition 4.19]{DGNO}.
\end{proof}

\subsection{Obstruction to fermionic actions}

In this subsection, we describe necessary and sufficient cohomological conditions for a group homomorphism $\rho\colon G \to \operatorname{Aut}_{\otimes}(\cC, f)$ to arise from a fermionic action of a super-group $(G, \alpha)$.

Let $\rho\colon G \to \operatorname{Aut}_{\otimes}(\cC, f)$ be a group homomorphism, where $\operatorname{Aut}_{\otimes}(\cC, f)$ denotes the group of isomorphism classes of tensor autoequivalences $F\colon \cC \to \cC$ such that $F(f) \cong f$ and diagram~\eqref{diagram:fermion functor} commutes. By Proposition \ref{prop:fermionic-functor-properties}, the commutativity of this diagram depends only on the underlying tensor functor $F$, not on the choice of isomorphism $F(f) \cong f$. In particular, every fermionic tensor functor has its underlying tensor functor in $\operatorname{Aut}_{\otimes}(\cC, f)$.

If there exists a fermionic action $\widetilde{\rho}_{\alpha}$ of $(G, \alpha)$ that induces $\rho$, we call $\widetilde{\rho}_{\alpha}$ an \emph{$\alpha$-lifting} of $\rho$.

In particular, $\widetilde{\rho}_{\alpha}$ defines a categorical action \smash{$\underline{G} \xrightarrow{\widetilde{\rho}_{\alpha}} \underline{\operatorname{Aut}_{\otimes}(\cC, f)} \subseteq \underline{\operatorname{Aut}_{\otimes}(\cC)}$} of $G$ that also realizes $\rho$. Hence, by Proposition \ref{prop:3-cocycle-obstruction}, the primary obstruction $O_3(\rho)$ must vanish.

Now consider a categorical action $\widetilde{\rho}\colon \underline{G} \to \underline{\operatorname{Aut}_{\otimes}(\cC, f)}$ of $G$ that realizes $\rho$. Since each $g_*(f)\cong f$, by Proposition \ref{diagram:fermion functor} there exist isomorphisms $\omega_g\colon g_*(f) \to f$ such that the monoidal structure satisfies \smash{$\psi(g)_{f,f}^{-1} = \omega_g \otimes \omega_g\colon g_*(f) \otimes g_*(f) \to \mathbf{1}$}.

We define the 2-cocycle $\theta_{\widetilde{\rho},\omega} \in Z^2(G, \mathbb{Z}/2\mathbb{Z})$ by requiring the commutativity of the following diagram:
\begin{equation}\label{diagrama theta_g}
\begin{tikzcd}
(gh)_*(f) \ar{rr}{\omega_{gh}} \ar{dd}[swap]{\phi(g,h)_{f}} && f \ar[ dotted]{dd}{(-1)^{\theta_{\widetilde{\rho},\omega}(g,h)}\id_f} \\\\
g_*h_*(f) \ar{rd}[swap]{g_*(\omega_h)} && f \\
& g_*(f) \ar{ur}[swap]{\omega_g}. &
\end{tikzcd}
\end{equation}

The diagram is well-defined since each $\omega_g$ and $\phi(g,h)_f$ are tensor natural isomorphisms. The dashed arrow represents the unique tensor automorphism of the identity functor that makes the diagram commutative, then must be multiplication by $\pm 1$.

\begin{Proposition}\label{prop: primera obstruccion}
Let \smash{$\widetilde{\rho}\colon \underline{G} \to \underline{\operatorname{Aut}_{\otimes}(\cC, f)}$} be a categorical action. Then
\begin{enumerate}\itemsep=0pt
 \item[$(1)$] The cohomology class of $\theta_{\widetilde{\rho},\omega} \in Z^2(G, \mathbb{Z}/2\mathbb{Z})$ does not depend on the choice of $\{\omega_g\}$. We denote this class by $[\theta_{\widetilde{\rho}}]$.
 \item[$(2)$] The action $\widetilde{\rho}$ is an $\alpha$-lifting if and only if $[\alpha] = [\theta_{\widetilde{\rho}}] \in H^2(G, \mathbb{Z}/2\mathbb{Z})$.
 \item[$(3)$] The fermionic structures on $\widetilde{\rho}$ form a torsor over $\operatorname{Hom}(G, \mathbb{Z}/2\mathbb{Z})$.
\end{enumerate}
\end{Proposition}

\begin{proof}
(1) Any other choice $\{\omega'_g\}$ differs by a normalized map $b\colon G \to \mathbb{Z}/2\mathbb{Z}$ via $\omega'_g\! =\!\smash{ (-1)^{b(g)}\omega_g}$, this $b$ relates $\delta(b)=\theta_{\widetilde{\rho},\omega}-\theta_{\widetilde{\rho},\omega'}$.

(2) A fermionic structure corresponds to a choice $\{u\colon G \to \mathbb{Z}/2\mathbb{Z}\}$ such that $\theta_{\widetilde{\rho},\omega}(g,h) - \alpha(g,h) = u(gh) - u(g) - u(h)$. Taking $\eta_g = (-1)^{u(g)}\omega_g$ gives the required fermionic structure.

Part (3) follows directly from the construction in part (2).
\end{proof}

The next step is to define the cohomological obstruction for the existence of $\alpha$-liftings in a~way that only depends on the data of $\rho$ and $\alpha$, and not on any particular choice of categorical lifting $\widetilde{\rho}$.

It follows from Proposition \ref{Proposition:Obstruction-bosonic} that the set of equivalence classes of liftings of $\rho$ is a torsor over \smash{$H^2_\rho\bigl(G, \widehat{K_0(\cC)}\bigr)$}. This means that if $\widetilde{\rho}$ is a fixed lifting, any other lifting $\widetilde{\rho}'$ can be written as $\widetilde{\rho}' = \beta \triangleright \widetilde{\rho}$ for some \smash{$\beta \in Z^2_\rho\bigl(G, \widehat{K_0(\cC)}\bigr)$}, where $\beta \triangleright \widetilde{\rho}$ means that we are twisting the natural isomorphism $\phi(g,h)\colon (gh)_*\to g_*\circ h_*$ as discussed in \eqref{eq: torciendo accion}.

By the definition of $\beta \triangleright \phi$ and diagram \eqref{diagrama theta_g}, we have
\smash{$
\theta_{\beta \triangleright \widetilde{\rho},\omega} = r_*(\beta) + \theta_{\widetilde{\rho},\omega}$},
where
\[
r_*\colon \ Z^2_\rho\bigl(G, \widehat{K_0(\cC)}\bigr) \to Z^2(G, \mathbb{Z}/2\mathbb{Z})
\]
 is the map induced by the restriction homomorphism \smash{$r\colon \widehat{K_0(\cC)} \to \mathbb{Z}/2\mathbb{Z}$} given by $(-1)^{r(\gamma)} = \gamma(f)$.

Hence, by Proposition \ref{prop: primera obstruccion} an $\alpha$-lifting exists if and only if for a given lifting $\widetilde{\rho}$, the equation
\begin{gather}\label{eq: r^*}
r_*(\beta) = [\alpha] - [\theta_{\widetilde{\rho}}]
\end{gather}
has a solution for some \smash{$\beta \in H^2_\rho\bigl(G, \widehat{K_0(\cC)}\bigr)$}.
The analysis of this equation splits into two cases.

\textit{Case $1$: The map $r$ is trivial.} In this case, the induced map $r_*$ is also trivial. The equation~${r_*(\beta) = [\alpha] + [\theta_{\widetilde{\rho}}]^{-1}}$ simplifies to $0 = [\alpha] + [\theta_{\widetilde{\rho}}]^{-1}$, which means a solution exists if and only if $[\alpha] = [\theta_{\widetilde{\rho}}]$.

\textit{Case $2$: The map $r$ is non-trivial.} To determine when a solution exists, we consider the short exact sequence of $G$-modules
\[
\xymatrix{
	1 \ar[r] & \operatorname{Ker}(r) \ar@{^{(}->}[r]^i & \widehat{K_0(\cC)} \ar@{->>}[r]^r & \mathbb{Z}/2\mathbb{Z} \ar[r] & 1
}
\]
and its associated long exact sequence in cohomology
\begin{gather*}
	\xymatrix{\ar[r] & H^2(G, \operatorname{Ker}(r)) \ar[r]^{i_*} & H^2\bigl(G,\widehat{K_0(\cC)}\bigr) \ar[r]^{r_*} & H^2(G, \mathbb{Z}/2\mathbb{Z}) \ar[r]^-{d_2} &} \\
	\qquad \xymatrix{\ar[r] & H^3(G, \operatorname{Ker}(r)) \ar[r]^{i_*} & H^3\bigl(G,\widehat{K_0(\cC)}\bigr) \ar[r]^{r_*} & H^3(G, \mathbb{Z}/2\mathbb{Z}) \ar[r]^-{d_3} &} \cdots.
\end{gather*}
By exactness, a solution for $\beta$ exists if and only if the element $[\alpha] - [\theta_{\widetilde{\rho}}]$ lies in the image of $r_*$, which is equivalent to requiring that it lies in the kernel of the connecting homomorphism $d_2$. This means $d_2([\alpha] - [\theta_{\widetilde{\rho}}]) = 0$.

The above analysis motivates the following definition of an obstruction class. We denote by~\smash{$\mathcal{K}_\rho := \operatorname{Ker}\bigl(r_*\colon H^2\bigl(G,\widehat{K_0(\cC)}\bigr) \to H^2(G,\mathbb{Z}/2\mathbb{Z})\bigr)$} the kernel of the induced map on cohomology.

\begin{Definition}\label{O3alt}
Let $\rho\colon G\to \operatorname{Aut}_\otimes(\cC,f)$ be a group homomorphism. The \emph{obstruction class} to the existence of an $\alpha$-lifting of $\rho$, denoted $O_3(\rho, \alpha)$, is defined as follows. Choose any lifting~${\widetilde{\rho}\colon \underline{G} \to \underline{\operatorname{Aut}_\otimes(\cC,f)}}$ of $\rho$ and set
\begin{enumerate}\itemsep=0pt
 \item[(i)] If the restriction map \smash{$r\colon \widehat{K_0(\cC)} \to \mathbb{Z}/2\mathbb{Z}$} is trivial, then $O_3(\rho, \alpha) := [\alpha] - [\theta_{\widetilde{\rho}}] \in H^2(G,\allowbreak\mathbb{Z}/2\mathbb{Z})$.
 \item[(ii)] If $r$ is non-trivial, then $O_3(\rho, \alpha) := d_2([\alpha] - [\theta_{\widetilde{\rho}}]) \in H^3(G,\operatorname{Ker}(r))$.
\end{enumerate}
This definition is independent of the choice of lifting $\widetilde{\rho}$ (see Theorem~\ref{theorem:bstruction fermionicaction}).
\end{Definition}

To state our next result, we require an additional construction. For each element ${a \in \mathcal{K}_\rho}$, we can choose a representative \smash{$\omega_a \in Z^2\bigl(G,\widehat{K_0(\cC)}\bigr)$} in the cohomology class $a$ and a 1-cochain $\gamma_a \in C^1(G,\mathbb{Z}/2\mathbb{Z})$ such that $r_*(\omega_a) = \delta(\gamma_a)$. This data determines a symmetric 2-cocycle
\begin{equation*}
\beta\colon\ \mathcal{K}_\rho \times \mathcal{K}_\rho \to H^1(G,\mathbb{Z}/2\mathbb{Z}),\qquad
(a,b) \mapsto \gamma_a + \gamma_b - \gamma_{a+b},
\end{equation*}
whose cohomology class in $H^2\bigl(\mathcal{K}_\rho, H^1(G,\mathbb{Z}/2\mathbb{Z})\bigr)$ is independent of the choice of representatives. This class defines an abelian group extension
\begin{equation}\label{eq:extension}
0 \to H^1(G,\mathbb{Z}/2\mathbb{Z}) \to \mathcal{A}_\rho \to \mathcal{K}_\rho \to 0.
\end{equation}

\begin{Theorem}\label{theorem:bstruction fermionicaction}
Let $(G,\alpha)$ be a super-group and $\rho\colon G \to \operatorname{Aut}_\otimes(\cC,f)$ a group homomorphism. The obstruction class $O_3(\rho, \alpha)$ is well-defined. Furthermore,
\begin{enumerate}[leftmargin=*,label={\rm (\alph*)}]\itemsep=0pt
 \item The homomorphism $\rho$ admits an $\alpha$-lifting if and only if $O_3(\rho,\alpha) = 0$.
 \item When $O_3(\rho,\alpha) = 0$, the set of equivalence classes of $\alpha$-liftings of $\rho$ $($as fermionic actions$)$ is a torsor over the abelian group $\mathcal{A}_\rho$ defined in \eqref{eq:extension}.
\end{enumerate}
\end{Theorem}

\begin{proof}
We first prove that $O_3(\rho, [\alpha])$ is well-defined, i.e., independent of the choice of lifting $\widetilde{\rho}$.

By Proposition \ref{Proposition:Obstruction-bosonic}, we know that any other lifting has the form $\widetilde{\rho}' = \beta \triangleright \widetilde{\rho}$ for some $[\beta] \in H^2_\rho\smash{\bigl(G, \widehat{K_0(\cC)}\bigr)}$, where
$
[\theta_{\widetilde{\rho}'}] = [\theta_{\beta \triangleright \widetilde{\rho}}] = r_*([\beta]) + [\theta_{\widetilde{\rho}}]$.

\textit{Case $1$: $r$ is trivial.} If $r$ is trivial, then $r_*$ is also trivial. It follows that $[\theta_{\widetilde{\rho}'}] = [\theta_{\widetilde{\rho}}]$ for any choice of $[\beta]$. Therefore, $O_3(\rho, [\alpha]) = [\alpha] - [\theta_{\widetilde{\rho}}]$ is independent of the choice of lifting.

\textit{Case $2$: $r$ is non-trivial.} We need to show that $d_2([\alpha] - [\theta_{\widetilde{\rho}'}]) = d_2([\alpha] - [\theta_{\widetilde{\rho}}])$. We compute:
\begin{align*}
d_2([\alpha] - [\theta_{\widetilde{\rho}'}]) &= d_2([\alpha] - (r_*([\beta]) + [\theta_{\widetilde{\rho}}])) = d_2(([\alpha] - [\theta_{\widetilde{\rho}}]) - r_*([\beta])) \\
&= d_2([\alpha] - [\theta_{\widetilde{\rho}}]) - d_2(r_*([\beta]))= d_2([\alpha] - [\theta_{\widetilde{\rho}}]).
\end{align*}
The last equality follows since $d_2 \circ r_* = 0$ by exactness of the long exact sequence. This shows that $O_3(\rho, [\alpha])$ is independent of the choice of lifting.

Next, we establish the characterization of when $\alpha$-liftings exist. From the preceding analysis, an $\alpha$-lifting exists if and only if equation \eqref{eq: r^*} has a solution
$ r_*([\beta]) = [\alpha] - [\theta_{\widetilde{\rho}}]$.

\textit{Case $1$: $r$ is trivial.} The equation becomes $0 = [\alpha] - [\theta_{\widetilde{\rho}}]$, which has a solution if and only if $[\alpha] = [\theta_{\widetilde{\rho}}]$. This is equivalent to $O_3(\rho, [\alpha]) = [\alpha] - [\theta_{\widetilde{\rho}}] = 0$.

\textit{Case $2$: $r$ is non-trivial.} A solution exists if and only if $[\alpha] - [\theta_{\widetilde{\rho}}] \in \operatorname{Im}(r_*)$. By exactness of the long exact sequence, this occurs if and only if $[\alpha] - [\theta_{\widetilde{\rho}}] \in \operatorname{Ker}(d_2)$, which means~${d_2([\alpha] - [\theta_{\widetilde{\rho}}]) = 0}$. This is precisely the condition $O_3(\rho, [\alpha]) = 0$.

When $O_3(\rho, [\alpha]) = 0$, there exists \smash{$[\beta] \in H^2_\rho\bigl(G,\widehat{K_0(\cC)}\bigr)$} such that $r_*([\beta]) = [\alpha] - [\theta_{\widetilde{\rho}}]$. Define~${\widetilde{\rho}' := (-[\beta]) \triangleright \widetilde{\rho}}$. Then
\begin{align*}
[\theta_{\widetilde{\rho}'}] &= r_*(-[\beta]) + [\theta_{\widetilde{\rho}}] = -r_*([\beta]) + [\theta_{\widetilde{\rho}}] = -([\alpha] - [\theta_{\widetilde{\rho}}]) + [\theta_{\widetilde{\rho}}] = [\alpha].
\end{align*}
Thus, $\widetilde{\rho}'$ is an $\alpha$-lifting of $\rho$.

Finally, we characterize the set of equivalence classes of $\alpha$-liftings. By Proposition \ref{prop: primera obstruccion}, fermionic structures on a fixed lifting of $\rho$ form a torsor over $H^1(G, \mathbb{Z}/2\mathbb{Z})$. Hence, it suffices to determine the different equivalence classes of $\alpha$-liftings as ordinary liftings of $\rho$ .

Suppose $\widetilde{\rho}_1$ and $\widetilde{\rho}_2$ are two $\alpha$-liftings of $\rho$. Then $\widetilde{\rho}_2 = [\beta] \triangleright \widetilde{\rho}_1$ for some \smash{$[\beta] \in H^2_\rho\bigl(G,\widehat{K_0(\cC)}\bigr)$}, and both satisfy $[\theta_{\widetilde{\rho}_1}] = [\theta_{\widetilde{\rho}_2}] = [\alpha]$. Using the relationship between different liftings
\begin{align*}
[\alpha] = [\theta_{\widetilde{\rho}_2}] = r_*([\beta]) + [\theta_{\widetilde{\rho}_1}] = r_*([\beta]) + [\alpha],
\end{align*}
which implies $r_*([\beta]) = 0$, hence $[\beta] \in \operatorname{Ker}(r_*)$.

Conversely, for any $[\beta] \in \operatorname{Ker}(r_*)$ and any $\alpha$-lifting $\widetilde{\rho}_\alpha$, the lifting $\widetilde{\rho}' = [\beta] \triangleright \widetilde{\rho}_\alpha$ satisfies $[\theta_{\widetilde{\rho}'}] = r_*([\beta]) + [\theta_{\widetilde{\rho}_\alpha}] = 0 + [\alpha] = [\alpha]$, showing that $\widetilde{\rho}'$ is also an $\alpha$-lifting.
\end{proof}

For a cyclic group $G = \mathbb{Z}/n\mathbb{Z}$, any $G$-module $M$, and $k > 0$, the cohomology groups $H^k_\rho(G, M)$ are periodic in $k$ with period 2. In particular, there is a canonical isomorphism
$
\smash{H^2_\rho(G, M)} \cong\smash{ \frac{M^G}{N(M)}}$,
where $M^G$ denotes the subgroup of elements fixed by the $\rho$-action of $G$, and $N\colon M \to M^G$ is the norm map defined by $
N(m) = \sum_{i=0}^{n-1} g_0^i \cdot m$
 for a generator $g_0$ of $G$.

For our application, the domain of $r_*$ is \smash{$H^2_\rho\bigl(G, \widehat{K_0(\cC)}\bigr) \cong \bigl(\widehat{K_0(\cC)}\bigr)^G / N\bigl(\widehat{K_0(\cC)}\bigr)$}. The codomain is $H^2(G, \mathbb{Z}/2\mathbb{Z})$ with trivial $G$-action. When $n$ is even, the norm map $N(x) = n \cdot x$ annihilates every element of $\mathbb{Z}/2\mathbb{Z}$, giving us $H^2(G, \mathbb{Z}/2\mathbb{Z}) \cong \mathbb{Z}/2\mathbb{Z}$.

\begin{Corollary}\label{cor:cyclic-fermionic-action}
Let $(\cC, f)$ be a fermionic fusion category and $G = \mathbb{Z}/n\mathbb{Z}$ with $n$ even. Consider a group homomorphism $\rho\colon G \to \operatorname{Aut}_\otimes(\cC, f)$ with vanishing obstruction $O_3(\rho) = 0 \in H^3_\rho\smash{\bigl(G, \widehat{K_0(\cC)}}\bigr)$.

Let $[\alpha] = 1 \in H^2(G, \mathbb{Z}/2\mathbb{Z}) \cong \mathbb{Z}/2\mathbb{Z}$ be the cohomology class associated to the non-trivial super-group $(\mathbb{Z}/2n\mathbb{Z}, [n])$. The existence of an $\alpha$-lifting of $\rho$ is characterized as follows:
\begin{enumerate}\itemsep=0pt
 \item[$(1)$] For any categorical lifting $\widetilde{\rho}$ of $\rho$, an $\alpha$-lifting exists if and only if the associated class $[\theta_{\widetilde{\rho}}] \in H^2(G, \mathbb{Z}/2\mathbb{Z})$ satisfies
 \[ 1 - [\theta_{\widetilde{\rho}}] \in \operatorname{Im}\left(r_*\colon \frac{\bigl(\widehat{K_0(\cC)}\bigr)^G}{N\bigl(\widehat{K_0(\cC)}\bigr)} \longrightarrow \mathbb{Z}/2\mathbb{Z}\right). \]
 This condition is independent of the choice of categorical lifting $\widetilde{\rho}$.
 \item[$(2)$] If $\rho$ is trivial, an $\alpha$-lifting exists if and only if $r_*$ is surjective.
 \item[$(3)$] When the existence criterion holds, the set of equivalence classes of $[\alpha]$-liftings is a torsor over the kernel $\mathcal{K} = \operatorname{Ker}(r_*)$.
\end{enumerate}
\end{Corollary}

\begin{proof}
(1) The general condition from Theorem~\ref{theorem:bstruction fermionicaction} for the existence of an $\alpha$-lifting is that $[\alpha] - [\theta_{\widetilde{\rho}}]$ must be in the image of $r_*$. For the non-trivial super-group, we have $[\alpha] = 1 \in \mathbb{Z}/2\mathbb{Z}$, and the condition becomes precisely the one stated in part (1).

(2) If $\rho$ is trivial, we can choose a categorical lifting $\widetilde{\rho}$ with $[\theta_{\widetilde{\rho}}] = 0$, and the conclusion follows from part (1).

(3) This follows directly from Theorem~\ref{theorem:bstruction fermionicaction}\,(b), which provides the classification of $\alpha$-liftings as a torsor over $\mathcal{K} = \operatorname{Ker}(r_*)$.
\end{proof}

\begin{Example}\label{ex:fermionic_actions_cyclic_trivial}
Let $\bigl(\operatorname{Vec}_G^\omega, (f, \eta)\bigr)$ be a fermionic pointed fusion category as in Example \ref{fermion pointed fusion categories}, where $G$ is a finite group and $f \in G$ is a central element of order 2. Let $n$ be even, and consider the non-trivial super-group $(\mathbb{Z}/n\mathbb{Z}, [\alpha])$ with $[\alpha] = 1 \in H^2(\mathbb{Z}/n\mathbb{Z}, \mathbb{Z}/2\mathbb{Z}) \cong \mathbb{Z}/2\mathbb{Z}$.

We apply Corollary~\ref{cor:cyclic-fermionic-action} to analyze the existence of $\alpha$-liftings that realize the trivial group homomorphism $\rho\colon \mathbb{Z}/n\mathbb{Z} \to \operatorname{Aut}_\otimes\bigl(\operatorname{Vec}_G^\omega, f\bigr)$.

By Corollary~\ref{cor:cyclic-fermionic-action}, an $\alpha$-lifting exists if and only if $r_*$ is surjective, where
\begin{align*}
r_*\colon\ \frac{\widehat{G}}{n\widehat{G}} \longrightarrow \mathbb{Z}/2\mathbb{Z}, \qquad
[\gamma] \mapsto \gamma(f).
\end{align*}
That is, an $[\alpha]$-lifting exists if and only if there exists a character $\gamma \in \widehat{G}$ such that ${\gamma(f) = -1}$. When such a $\alpha$-lifting exists, the set of equivalence classes is a torsor over $\mathcal{K}\oplus \mathbb{Z}/2\mathbb{Z}$ where \smash{$\mathcal{K} = \operatorname{Ker}(r_*) = \frac{f^\perp }{n\widehat{G}}$},
and $f^\perp := \{\gamma \in \widehat{G} \mid \gamma(f) = 1\}$.
\end{Example}

A super-group of the form $\widetilde{G}=G \times \mathbb{Z}/2\mathbb{Z}$ is called a \emph{trivial super-group}. Equivalently, a~trivial super-group corresponds to a pair $(G,\alpha)$ with $[\alpha]=0 \in H^2(G,\mathbb{Z}/2\mathbb{Z})$ the trivial class.

\begin{Corollary}\label{coro:trivial-action}
Let $(\cC,f)$ be a fermionic fusion category and $G$ a finite group. The trivial group homomorphism $\operatorname{triv}\colon G\to \operatorname{Aut}_\otimes(\cC,f)$ extends to a fermionic action of a super-group~${(G,\alpha)}$ if and only if
\begin{enumerate}[leftmargin=*,label={\rm (\alph*)}]\itemsep=0pt
 \item If the restriction map \smash{$r\colon \widehat{K_0(\cC)} \to \mathbb{Z}/2\mathbb{Z}$} is trivial, then $\alpha=0$.
 \item If $r$ is non-trivial, then $\alpha \in \operatorname{Im}(r_*)$.
 \item If the epimorphism $r$ from part (b) splits as a map of $G$-modules, then a fermionic action exists for every super-group $(G,\alpha)$.
\end{enumerate}
\end{Corollary}

\begin{proof}
By Theorem~\ref{theorem:bstruction fermionicaction}, this occurs if and only if $O_3(\operatorname{triv}, \alpha) = 0$. For the trivial homomorphism, we can choose a categorical lifting with $[\theta_{\widetilde{\operatorname{triv}}}] = 0$.

(a) If $r$ is trivial, then $O_3(\operatorname{triv}, \alpha) = \alpha$, so the condition becomes $\alpha = 0$.

(b) If $r$ is non-trivial, then $O_3(\operatorname{triv}, \alpha) = d_2(\alpha)$. By exactness of the long exact sequence in cohomology, $d_2(\alpha) = 0$ if and only if $\alpha \in \operatorname{Im}(r_*)$.

(c) If \smash{$r\colon \widehat{K_0(\cC)} \to \mathbb{Z}/2\mathbb{Z}$} splits, then the induced map \[
r_*\colon\ H^2\big(G, \widehat{K_0(\cC)}\big) \to H^2(G, \mathbb{Z}/2\mathbb{Z})
\]
 is surjective. Therefore, $\operatorname{Im}(r_*) = H^2(G, \mathbb{Z}/2\mathbb{Z})$, and the condition $\alpha \in \operatorname{Im}(r_*)$ holds for any $\alpha \in H^2(G, \mathbb{Z}/2\mathbb{Z})$.
\end{proof}

\begin{Example}
\label{ex:elementary-abelian-2-groups}
Let $A = (\mathbb{Z}/2\mathbb{Z})^n$ and \smash{$\bigl(\operatorname{Vec}_A^\omega, (f, \eta)\bigr)$} a fermionic pointed fusion category with~${f \in A}$ of order~2. The restriction map $r\colon \widehat{A} \to \mathbb{Z}/2\mathbb{Z}$ given by $(-1)^{r(\gamma)} = \gamma(f)$ splits.

By Corollary~\ref{coro:trivial-action}\,(c), every super-group $(G, \alpha)$ admits a fermionic action on $\bigl(\operatorname{Vec}_A^\omega, f\bigr)$ via the trivial homomorphism.
\end{Example}

\subsection{Group-theoretical interpretation of the obstruction}

The goal of this subsection is to provide a group-theoretical interpretation of the obstruction~${O_3(\rho,\alpha)}$ from Theorem~\ref{theorem:bstruction fermionicaction}. While the cohomological criterion $O_3(\rho,\alpha) = 0$ is clear, its verification in concrete examples often requires understanding the underlying group-theoretic structure. We translate the cohomological condition into a statement about morphisms between group extensions, which allows us to exploit structural properties of specific groups (such as abelianness or centrality conditions) to determine when fermionic actions exist.

It is a classical result, due to Schreier and Eilenberg--Mac~Lane \cite{EM1}, that extensions of a~group~$G$ by a $G$-module $C$ are classified by the second cohomology group $H^2(G,C)$. Specifically, there is a canonical bijection between elements of $H^2(G,C)$ and equivalence classes of group extensions
$
1 \longrightarrow C \longrightarrow E \longrightarrow G \longrightarrow 1$,
where $E$ is a group and the action of $G$ on $C$ induced by conjugation in $E$ coincides with the given $G$-module structure. Two extensions are equivalent if there exists a group isomorphism~${E \to E'}$ making the obvious diagram commute and inducing the identity on both $C$ and $G$.

The problem of lifting a cohomology class via a module epimorphism can now be stated in terms of these group extensions.

\begin{Lemma}\label{lemma: diagrama para poder extender}
Let $r\colon B \to C$ be an epimorphism of $G$-modules and let $[\alpha] \in H^2(G,C)$. A~class $[\beta] \in H^2(G,B)$ satisfying $r_*([\beta])=[\alpha]$ exists if and only if there is a group epimorphism~${\widetilde{r}\colon G_\beta \to G_\alpha}$ such that the following diagram commutes
\begin{equation}\label{cuadro:cuadro para existencia de imagen de restriccion}
\vcenter{\hbox{\xymatrix{
 0 \ar[r] & B \ar@{^{(}->}[r] \ar@{->>}[d]_{r} & G_\beta \ar@{->>}[r] \ar@{->>}[d]^{\widetilde{r}} & G \ar@{=}[d] \ar[r] & 0 \\
 0 \ar[r] & C \ar@{^{(}->}[r] & G_\alpha \ar@{->>}[r] & G \ar[r] & 0.
}}}
\end{equation}
\end{Lemma}

Consider the case where $(G,\alpha)$ is a super-group and $\rho\colon G \to \operatorname{Aut}_\otimes(\cC,f)$ is a group homomorphism with $r$ non-trivial. Recall from the proof of Theorem~\ref{theorem:bstruction fermionicaction} that the existence of an $\alpha$-lifting depends on the relationship between $[\alpha]$ and $[\theta_{\widetilde{\rho}}]$ via the connecting homomorphism $d_2$. We now interpret this condition in terms of group extensions.

The short exact sequence of $G$-modules
$
\xymatrix{
 0 \ar[r] & \operatorname{Ker}(r) \ar[r]^{i} & \widehat{K_0(\cC)} \ar@{->>}[r]^r & \mathbb{Z}/2\mathbb{Z} \ar[r] & 0
}
$
induces a long exact sequence in group cohomology. By Lemma \ref{lemma: diagrama para poder extender}, the cohomological condition~${O_3(\rho,\alpha) = 0}$ is equivalent to the existence of a commutative diagram of group extensions. Specifically, we seek an extension of $G$ by \smash{$\widehat{K_0(\cC)}$ }such that the diagram
\[
\xymatrix{
 0 \ar[r] & \widehat{K_0(\cC)} \ar@{^{(}->}[r] \ar@{->>}[d]^{r} & L \ar@{->>}[r] \ar@{->>}[d]^{\widetilde{r}} & G \ar@{=}[d] \ar[r] & 0 \\
 0 \ar[r] & \mathbb{Z}/2\mathbb{Z}\ar@{^{(}->}[r] & G_\varphi \ar@{->>}[r]^{\pi_\varphi} & G \ar[r] & 0
}
\]
commutes. Here, $G_\varphi$ is the super-group corresponding to $\varphi = \theta_{\widetilde{\rho}}/\alpha \in H^2(G,\mathbb{Z}/2\mathbb{Z})$.

The following proposition provides additional structure when $G$ acts trivially on $\operatorname{Ker}(r)$.

\begin{Proposition}\label{proposition: central extensions for extensions}
Consider the commutative diagram \eqref{cuadro:cuadro para existencia de imagen de restriccion} of Lemma {\rm\ref{lemma: diagrama para poder extender}}. If $G$ acts trivially on $\operatorname{Ker}(r)$, then
\begin{enumerate}\itemsep=0pt
 \item[$(1)$] $\operatorname{Ker}(r)$ is a central subgroup of $G_\beta$.
 \item[$(2)$] The following diagram is commutative{\samepage
 \[
\begin{tikzcd}
	& 0 \arrow[d] & 0 \arrow[d] \\
	& \operatorname{Ker}(r) \arrow[r, equals] \arrow[d] & \operatorname{Ker}(r) \arrow[d] \\
	0 \arrow[r] & B \arrow[r, "i_{G_{\beta}}", hook] \arrow[d, "r"] & G_{\beta} \arrow[r, "\pi_{\beta}"] \arrow[d, "\tilde{r}"] & G \arrow[d, equals] \arrow[r] & 0 \\
	0 \arrow[r] & C \arrow[r, "i_{G_{\alpha}}", hook] \arrow[d] & G_{\alpha} \arrow[r, "\pi_{\alpha}"'] \arrow[d] & G \arrow[r] & 0, \\
	& 0 & 0 & &
\end{tikzcd}
 \]
where both columns are central extensions.}
 \item[$(3)$] If $B$ corresponds to $\delta \in H^2(C,\operatorname{Ker}(r))$ and $G_\beta$ to $\psi \in H^2(G_\alpha,\operatorname{Ker}(r))$, then $\delta = i_{G_\alpha}^*(\psi)$.
\end{enumerate}
\end{Proposition}

\begin{Remark}
Proposition \ref{proposition: central extensions for extensions} reduces the existence problem for fermionic actions to purely group-theoretic obstructions. When $G$ acts trivially on $\operatorname{Ker}(r)$, the diagram \eqref{cuadro:cuadro para existencia de imagen de restriccion} forces $G_\beta$ to have a central subgroup of a specific form. For non-abelian super-groups or non-trivial actions of $G$ on \smash{$\widehat{K_0(\cC)}$}, these centrality requirements often lead to contradictions, as we illustrate in the following examples.
\end{Remark}

We now apply these results to study specific examples of fermionic actions.

\begin{Example}\label{example: fermionic action non trivial z2n}
Consider $\operatorname{Vec}_A^q$ with fermion $f$, a pointed spin-braided category from Example \ref{ex: pointed 4 rank}, and the super-group $(\mathbb{Z}/2n\mathbb{Z},[n])$ with a non-trivial homomorphism $\rho\colon \mathbb{Z}/n\mathbb{Z} \to \operatorname{Aut}_\otimes\bigl(\operatorname{Vec}_A^q,f\bigr)$. By Theorem~\ref{theorem: extensions of A with order 4}, $O_3(\rho) = 0$.

Since $\mathbb{Z}/n\mathbb{Z}$ acts trivially on $\operatorname{Ker}(r) \cong \mathbb{Z}/2\mathbb{Z}$, by Lemma \ref{lemma: diagrama para poder extender} and Proposition \ref{proposition: central extensions for extensions}, $\rho$ has an $\alpha$-lifting if and only if there exists a commutative diagram
\[\label{diagram: super-actions of ciclic groups}
\begin{tikzcd}
	& 0 \arrow[d] & 0 \arrow[d] \\
	& \operatorname{Ker}(r) \arrow[r, equals] \arrow[d] & \operatorname{Ker}(r) \arrow[d] \\
	0 \arrow[r] & \widehat{A} \arrow[r, hook] \arrow[d, "r", two heads] & L \arrow[r, two heads] \arrow[d, "\widetilde{r}", two heads] & \mathbb{Z}/n\mathbb{Z} \arrow[d, equals] \arrow[r] & 0 \\
	0 \arrow[r] & \mathbb{Z}/2\mathbb{Z} \arrow[r, hook] \arrow[d] & \mathbb{Z}/2n\mathbb{Z} \arrow[r, "\pi_{\varphi}", two heads] \arrow[d] & \mathbb{Z}/n\mathbb{Z} \arrow[r] & 0, \\
	& 0 & 0 & &
\end{tikzcd}
\]
where the columns are central extensions. Since middle column is a central extension by cyclic group, $L$ must be abelian. However, since the action of $\mathbb{Z}/n\mathbb{Z}$ on $\widehat{A}$ is non-trivial, $L$ cannot be abelian, yielding a contradiction. Therefore, the super-group $(\mathbb{Z}/2n\mathbb{Z},[n])$ admits no fermionic action on $\operatorname{Vec}_A^q$ realizing the non-trivial homomorphism $\rho$.
\end{Example}

\begin{Example}
\label{ex:no-action-D4-on-Z4}
Let $\bigl(\widetilde{G}, z\bigr)$ be a non-abelian super-group of order 8, hence $\widetilde{G}$ is either the dihedral group $D_4$ or the quaternion group $Q_8$, and $z$ is the unique non-trivial central element of order 2. Then $G := \widetilde{G}/\langle z \rangle \cong \mathbb{Z}/2\mathbb{Z} \times \mathbb{Z}/2\mathbb{Z}$.

We show that this super-group does not admit a fermionic action on $\bigl(\operatorname{Vec}_{\mathbb{Z}/4\mathbb{Z}}^q, f\bigr)$ with fermion~${f = 2}$ via the trivial homomorphism.

By Lemma \ref{lemma: diagrama para poder extender}, such an action would require a group $L$ of order 16 fitting into the commutative diagram
\[
\begin{tikzcd}
 0 \arrow[r] & \mathbb{Z}/4\mathbb{Z} \arrow[r] \arrow[d, "r", two heads] & L \arrow[r, two heads] \arrow[d, "\tilde{r}", two heads] & G \arrow[d, equals] \arrow[r] & 0 \\
 0 \arrow[r] & \mathbb{Z}/2\mathbb{Z} \arrow[r] & \widetilde{G} \arrow[r, two heads] & G \arrow[r] & 0.
\end{tikzcd}
\]

From this diagram, $L$ must satisfy two contradictory properties:
\begin{enumerate}\itemsep=0pt
 \item[(1)] $L$ has a central subgroup $\mathbb{Z}/4\mathbb{Z} \subseteq Z(L)$ (since the top row is a central extension for trivial~$\rho$).
 \item[(2)] $L$ is non-abelian (since $\tilde{r}\colon L \to \widetilde{G}$ is surjective and $\widetilde{G}$ is non-abelian).
\end{enumerate}

However, no non-abelian group of order 16 has center containing $\mathbb{Z}/4\mathbb{Z}$. All such groups have center isomorphic to $\mathbb{Z}/2\mathbb{Z}$ or $\mathbb{Z}/2\mathbb{Z} \times \mathbb{Z}/2\mathbb{Z}$.

This contradiction shows that $\bigl(\widetilde{G}, z\bigr)$ admits no fermionic action on $\bigl(\operatorname{Vec}_{\mathbb{Z}/4\mathbb{Z}}^q, f\bigr)$.
\end{Example}

\begin{Example}[no non-trivial fermionic action of $(Q_8, z)$]
\label{ex:no-action-Q8-on-Z4}
Let $(Q_8, z)$ be the quaternion super-group of order 8. We show that it admits no fermionic action on $(\operatorname{Vec}_{A}^q, f)$, a non-degenerate spin-braided pointed fusion category of rank 4 as in Example \ref{ex: pointed 4 rank}, via a non-trivial homomorphism~${\rho\colon \mathbb{Z}/2\mathbb{Z} \times \mathbb{Z}/2\mathbb{Z} \to \operatorname{Aut}_\otimes\bigl(\operatorname{Vec}_{A}^q, f\bigr)}$.

We proceed by contradiction. Assume such an action exists. By Lemma \ref{lemma: diagrama para poder extender}, there would exist a group $L$ of order 16 with properties:
\begin{enumerate}\itemsep=0pt
 \item[(1)] $L$ contains a normal subgroup \smash{$\widehat{A}$} with non-trivial $G$-action.
 \item[(2)] The central subgroup $ \mathbb{Z}/2\mathbb{Z}\cong K := \operatorname{Ker}(r)\subset \widehat{A} \subset L
 $ satisfies $L/K \cong Q_8$.
\end{enumerate}

Let $S = \big\langle \widehat{A}, y \big\rangle$ where $y \in L$ does not centralize $\widehat{A}$. Then $S$ has order 8 and is not abelian, hence it is isomorphic to either $D_4$ or $Q_8$. Consider the projection $\pi\colon L \to L/K \cong Q_8$. Since $|S \cap K| = 2$, we have $|\pi(S)| = 4$.

The quotient $\pi(S) \cong S/K$ is isomorphic to $\mathbb{Z}/2\mathbb{Z} \times \mathbb{Z}/2\mathbb{Z}$, since both $D_4/Z(D_4)$ and $Q_8/Z(Q_8)$ equal the Klein four-group.

However, all subgroups of order 4 in $Q_8$ are cyclic. This contradiction shows that $(Q_8, z)$ admits no non-trivial fermionic action on $\bigl(\operatorname{Vec}_{A}^q, f\bigr)$.
\end{Example}

\subsection[Fermionic actions on non-degenerate spin-braided fusion categories of rank four]{Fermionic actions on non-degenerate spin-braided fusion categories\\ of rank four}\label{accionsv}

In this subsection, we apply the results developed earlier to study fermionic actions on non-degenerate spin-braided fusion categories of dimension four. We first characterize fermionic actions on Ising braided categories, then extend our analysis to pointed spin-braided fusion categories of rank four.

\subsubsection{Fermionic actions on Ising braided categories}

We begin by investigating which super-groups can act on Ising braided categories, see Example~\ref{Example:Ising}.

\begin{Lemma}\label{lemma:kernel r in spin-modular case}
Let $(\cB,f)$ be a non-degenerate spin-braided category. The restriction map $r\colon\allowbreak\smash{\widehat{K_0(\cB)}} \to \mathbb{Z}/2\mathbb{Z}$ is trivial if and only if $\operatorname{Inv}(\cB)=\operatorname{Inv}(C_{\cB}(f))$, where $C_{\cB}(f)$ is the centralizer of~$f$ in $\cB$, consisting of objects $X$ such that $c_{f,X}\circ c_{X,f}=\id_{X \otimes f}$.
\end{Lemma}

\begin{proof}
If $\cB$ is a non-degenerate braided fusion category, there is a canonical isomorphism \smash{$\Xi\colon \operatorname{Inv}(\cB) \to \widehat{K_0(\cB)}$} given by $\Xi(X_i)(X_j)=c_{X_j,X_i} \circ c_{X_i,X_j}$ (see \cite[Section 6.1]{GELAKI20081053}).

For a non-degenerate spin-braided category $(\cB,f)$, consider the map $s$ defined by the composition
\[
\xymatrix{
 \widehat{K_0(\cB)} \ar[r]^r & \mathbb{Z}/2\mathbb{Z} \\
 \operatorname{Inv}(\cB) \ar[u]^{\Xi} \ar@{..>}[ur]_s.
}
\]

The kernel of $s$ consists of all isomorphism classes of invertible objects $X$ such that $c_{X,f}\circ c_{f,X}=\id_{f \otimes X}$, which by definition is $\operatorname{Inv}(C_{\cB}(f))$. Thus, the restriction map $r$ is trivial if and only if $\operatorname{Inv}(\cB)=\operatorname{Inv}(C_{\cB}(f))$.
\end{proof}

Lemma \ref{lemma:kernel r in spin-modular case} provides a criterion for determining when the obstruction theory for fermionic actions simplifies. We now apply it to Ising categories.

\begin{Proposition}\label{proposition: fermionicactionsonIsing}
Only trivial super-groups $(G,\alpha \equiv 0)$ act fermionically on spin-braided Ising categories. Moreover, the equivalence classes of fermionic actions of a group $G$ on an Ising category are in bijective correspondence with elements of $H^2(G,\mathbb{Z}/2\mathbb{Z})$.
\end{Proposition}

\begin{proof}
Ising categories are particular cases of Tambara--Yamagami categories. By \cite[Proposition 1]{Tambara}, every tensor auto-equivalence of an Ising category $\mathcal{I}$ is monoi\-dally equivalent to the identity functor. Thus, $\operatorname{Aut}_\otimes(\mathcal{I},f) \cong \{\id_{\mathcal{I}}\}$.

Since $\operatorname{Inv}(\mathcal{I}_{\rm pt}) = \langle f \rangle$ for an Ising category $\mathcal{I}$, by Lemma \ref{lemma:kernel r in spin-modular case}, the homomorphism $r\colon \smash{\widehat{K_0(\mathcal{I})}} \to \mathbb{Z}/2\mathbb{Z}$ is trivial. Therefore, by Corollary~\ref{coro:trivial-action}\,(a), only trivial super-groups can act on $\mathcal{I}$.

For a trivial super-group $(G,0)$, the set of equivalence classes of fermionic actions is a torsor over \smash{$\operatorname{Ker}(r_*) = H^2\bigl(G,\widehat{K_0(\mathcal{I})}\bigr)$}. Since every tensor auto-equivalence of $\mathcal{I}$ is equivalent to the identity, we have \smash{$\widehat{K_0(\mathcal{I})} \cong \mathbb{Z}/2\mathbb{Z}$} with trivial $G$-action. Thus, the equivalence classes of fermionic actions are in bijective correspondence with $H^2(G,\mathbb{Z}/2\mathbb{Z})$.
\end{proof}

\subsubsection{Fermionic actions on pointed spin-braided fusion categories of rank four}\label{subsubsection:fermionic actions on pointed spin-modular categories}

We now turn our attention to pointed non-degenerate spin-braided fusion categories of dimension four described in Example \ref{ex: pointed 4 rank}.

\begin{Theorem}\label{theorem: extensions of A with order 4}
Let $(G,\alpha)$ be a finite super-group and $\bigl(\operatorname{Vec}_A^q,f\bigr)$ a pointed spin-braided fusion category of rank four as described in Example {\rm\ref{ex: pointed 4 rank}}. Then
\begin{enumerate}[leftmargin=*,label={\rm (\alph*)}]\itemsep=0pt
 \item $\operatorname{Aut}_\otimes^{\rm br}\bigl(\operatorname{Vec}_A^q,f\bigr) = O(A,q,f) \cong \mathbb{Z}/2\mathbb{Z}$.
 \item Any group homomorphism $\rho\colon G\to O(A,q,f)$ admits a categorical lifting
 \[
 \widetilde{\rho}\colon \smash{\underline{G} \to \underline{\operatorname{Aut}_\otimes^{\rm br}\bigl(\operatorname{Vec}_A^q,f\bigr)}}
 \]
 with $\theta_{\widetilde{\rho}} = 1$.
 \item $\rho$ can be realized by a fermionic action of $(G,\alpha)$ if and only if $d_2(\alpha)=0$.
 \item If $d_2(\alpha)=0$, the equivalence classes of such fermionic actions form a torsor over $\operatorname{Ker}\bigl(r_*\colon\allowbreak H^2\bigl(G,\widehat{A}\bigr) \to H^2(G, \mathbb{Z}/2\mathbb{Z})\bigr)$.
\end{enumerate}
Here $d_2\colon H^2(G,\mathbb{Z}/2\mathbb{Z})\to H^3(G,\mathbb{Z}/2\mathbb{Z})$ is the connecting homomorphism for the exact sequence $0\to f^\perp \to \widehat{A} \overset{r}{\to} \mathbb{Z}/2\mathbb{Z} \to 0$, where $f^\perp := \big\{\gamma \in \widehat{A} \mid \gamma(f) = 1\big\}$.
\end{Theorem}

\begin{proof}
(a) By Example \ref{ex: pointed 4 rank}, $O(A,q,f) \cong \mathbb{Z}/2\mathbb{Z}$ consists of automorphisms preserving both $q$ and $f$. By definition, $\operatorname{Aut}_\otimes^{\rm br}\bigl(\operatorname{Vec}_A^q,f\bigr) = O(A,q,f)$.

(b) We identify roots of unity in $\mathbb{C}$ with $\mathbb{Q}/\mathbb{Z}$. Following Examples \ref{fermionic action on pointed fusion categories} and \ref{ex:fermionic_actions_cyclic_trivial}, we construct explicit categorical actions. Given an abelian 3-cocycle $(\omega,c)$ realizing $q$, a braided action of~$\mathbb{Z}/2\mathbb{Z}$ on $\operatorname{Vec}_A^q$ is determined by maps $\mu\colon A \times A \to \mathbb{C}^\times$ and $\gamma\colon A \to \mathbb{C}^\times$ satisfying
\begin{gather*}
\omega(a,b,c) - \omega(g_*(a), g_*(b), g_*(c)) = \mu(b,c) + \mu(a,b+c) - \mu(a+b,c) - \mu(a,b), \\
\mu(g_*(a), g_*(b)) + \mu(a,b) = \gamma(ab) - \gamma(a) - \gamma(b), \qquad
\gamma(g_*(a)) = \gamma(a), \\
c(a,b)-c(g_*(a),g_*(b)) = \mu(a,b)-\mu(b,a),
\end{gather*}
where $g$ generates $\mathbb{Z}/2\mathbb{Z}$.

We exhibit the data $(\omega,c,\mu,\gamma)$ for each case in order:

\textit{Case $1$a:} $A = \mathbb{Z}/2\mathbb{Z} \times \mathbb{Z}/2\mathbb{Z}$, $q_0((x,y)) = \frac{xy}{2}$, $g(x,y) = (y,x)$.

Data: $c_0((x_1,x_2),(y_1,y_2)) = \frac{x_1y_2}{2}$, $\omega = 0$, $\mu_g((x_1,x_2),(y_1,y_2)) = \frac{x_1y_1 + x_2y_1}{2}$.
\begin{center}
\begin{tabular}{c|c}
$a$ & $\gamma_g(a)$ \\
\hline
$(0,0)$ & $0$ \\
$(1,0)$ & $\frac{1}{4}$ \\
$(0,1)$ & $\frac{1}{4}$ \\
$(1,1)$ & $0$
\end{tabular}
\end{center}

\textit{Case $1$b:} $A = \mathbb{Z}/2\mathbb{Z} \times \mathbb{Z}/2\mathbb{Z}$, $q_{4}((x,y)) = \frac{x^2+xy+y^2}{2}$, $g(x,y) = (y,x)$.

Data: $c_{1/2}((x_1,x_2),(y_1,y_2)) = \frac{x_1y_1+x_2y_2+x_1y_2}{2}$, $\omega = 0$. The maps $\mu_g$ and $\gamma_g$ are identical to Case 1a.

\textit{Case $1$c:} $A = \mathbb{Z}/2\mathbb{Z} \times \mathbb{Z}/2\mathbb{Z}$, $q_\pm((x,y)) = \pm\frac{x^2+y^2}{4}$, $g(x,y) = (y,x)$.

Data: $c((x_1,x_2),(y_1,y_2)) = \pm\bigl(\frac{x_1y_1}{4} + \frac{x_2y_2}{4}\bigr)$, \[
\omega((x_1,y_1),(x_2,y_2),(x_3,y_3)) = \frac{x_1}{2}\left\lfloor \frac{y_1+z_1}{4} \right\rfloor + \frac{x_2}{2}\left\lfloor \frac{y_2+z_2}{4} \right\rfloor,
\]
 $\mu_g \equiv 0$, $\gamma_g \equiv 0$.

\textit{Case $2$:} $A = \mathbb{Z}/4\mathbb{Z}$, $q_k(\bar{x}) = \frac{kx^2}{8}$, $g(\bar{x}) = -\bar{x}$ and $k\in \{1,3,5,7\}$.

Data: $c(x,y) = \frac{kxy}{8}$, $\omega(x,y,z) = \frac{x}{2}\big\lfloor \frac{y+z}{4} \big\rfloor$, $\gamma_g \equiv 0$.

\begin{center}
\begin{tabular}{c|cccc}
$\mu_g(\bar{a},\bar{b})$ & $\bar{0}$ & $\bar{1}$ & $\bar{2}$ & $\bar{3}$ \\
\hline
$\bar{0}$ & $0$ & $0$ & $0$ & $0$ \\
$\bar{1}$ & $0$ & $0$ & $\frac{1}{2}$ & $0$ \\
$\bar{2}$ & $0$ & $0$ & $0$ & $0$ \\
$\bar{3}$ & $0$ & $0$ & $\frac{1}{2}$ & $0$
\end{tabular}
\end{center}

For any $\rho\colon G \to O(A,q,f)$, the composition
\begin{equation*}
\begin{tikzcd}
\underline{G} \arrow[d, "\pi"'] \arrow[r, dashed, "\widetilde{\rho}"] & \underline{\operatorname{Aut}_\otimes^{\rm br}\bigl(\operatorname{Vec}_A^q,f\bigr)} \\
\underline{\mathbb{Z}/2\mathbb{Z}} \arrow[ur, "F"']
\end{tikzcd}
\end{equation*}
yields a categorical lifting with $\theta_{\widetilde{\rho}} = 1$.

(c) By Definition~\ref{O3alt}, since $\theta_{\widetilde{\rho}} = 1$ and the restriction map $r\colon \widehat{A} \to \mathbb{Z}/2\mathbb{Z}$ given by $r(\gamma) = \gamma(f)$ is non-trivial (as $|A| = 4$ and $f$ has order 2), we have $O_3(\rho,\alpha) = d_2(\alpha)$. Thus $\rho$ admits a~fermionic $(G,\alpha)$-action if and only if $d_2(\alpha) = 0$.

(d) When $d_2(\alpha) = 0$, by Theorem~\ref{theorem:bstruction fermionicaction}\,(b), the equivalence classes form a torsor over $\operatorname{Ker}\bigl(r_*\colon H^2\bigl(G,\widehat{A}\bigr) \to H^2(G, \mathbb{Z}/2\mathbb{Z})\bigr)$.
\end{proof}

This theorem completely characterizes which super-groups can act on non-degenerate spin-braided pointed fusion categories of rank four, and classifies the possible actions. The condition~${d_2(\alpha)=0}$ provides a cohomological criterion for the existence of such actions.

An important consequence is that non-trivial super-groups may act on pointed spin-braided fusion categories, in contrast to the Ising case where only trivial super-groups can act. This highlights a fundamental difference between pointed and non-pointed spin-braided fusion categories of dimension four.

\section{Minimal non-degenerate extensions of super-groups} \label{section: minimal non-degenerated extension}

In this section, we study minimal non-degenerate extensions of braided fusion categories, with special emphasis on super-groups. We begin by reviewing relevant concepts about module categories and the Brauer--Picard group, then establish necessary and sufficient conditions for the existence of minimal non-degenerate extensions.

\subsection{Module categories and Brauer--Picard groups}

Let $\cC$ be a fusion category. A left $\cC$-\emph{module category} is a semisimple $\mathbb{C}$-linear abelian category~$\cM$ equipped with a $\mathbb{C}$-bilinear bi-exact bifunctor $\otimes\colon \cC \times \cM \to \cM$ and natural associativity and unit isomorphisms $m_{X,Y,M}\colon (X \otimes Y) \otimes M \to X \otimes (Y \otimes M)$, $\lambda_M\colon \mathbf{1} \otimes M \to M$ satisfying the usual coherence conditions.
Right module categories and bimodule categories over $\cC$ are defined similarly. The tensor product $\boxtimes_{\cC}$ of $\cC$-bimodule categories was defined in \cite{ENO3}, giving the 2-category of $\cC$-bimodule categories the structure of a monoidal bicategory \cite{tricategories}.
A bimodule category $\cM$ is called \emph{invertible} if there exists a $\cC$-bimodule $\cN$ such that ${\cM \boxtimes_{\cC} \cN \cong \cC}$ and $\cN \boxtimes_{\cC} \cM \cong \cC$ as $\cC$-bimodule categories. The \emph{Brauer--Picard group} $\operatorname{BrPic}(\cC)$ of $\cC$ is the group of equivalence classes of invertible $\cC$-bimodule categories. By \cite[Theorem~1.3]{ENO3}, this group plays a key role in the classification of extensions of tensor categories by finite groups.

The natural structure for invertible bimodule categories over a fusion category $\cC$ is the $3$-group \smash{$\underline{\underline{\operatorname{BrPic}(\cC)}}$}, whose $1$-truncation is the $2$-group \smash{$\underline{\operatorname{BrPic}(\cC)}$}, and whose $2$-truncation is the Brauer--Picard group $\operatorname{BrPic}(\cC)$.

For a braided fusion category $\cB$, a left action induces a compatible right action via the braiding. Thus, all left~$\cB$-modules have a canonical $\cB$-bimodule structure. In the braided case, there is a distinguished $3$-subgroup \smash{$\underline{\underline{\operatorname{Pic}(\cB)}} \subseteq \underline{\underline{\operatorname{BrPic}(\cB)}}$}, called the \emph{Picard $3$-group} of $\cB$, consisting of all invertible (left) $\cB$-modules.

Let \smash{$\underline{\operatorname{Aut}_\otimes^{\rm br}(\cB)}$} be the 2-group of braided tensor auto-equivalences of a braided fusion category~$\cB$. There is a monoidal functor \smash{$\Theta\colon \underline{\operatorname{Pic}(\cB)} \to \underline{\operatorname{Aut}_\otimes^{\rm br}(\cB)}$} associated with the \emph{alpha-induction functors} $\alpha_{+}$ and $\alpha_{-}$ (see \cite{ENO3,ostrik1} for precise definitions). When $\cB$ is non-degenerate, there also exists a monoidal functor \smash{$\Phi\colon \underline{\operatorname{Aut}_\otimes^{\rm br}(\cB)} \to \underline{\operatorname{Pic}(\cB)}$} such that the functors $\Theta$ and $\Phi$ are mutually inverse monoidal equivalences \cite[Theorem 5.2]{ENO3}. Using this equivalence, we will identify \smash{$\underline{\operatorname{Pic}(\cB)}$} with \smash{$\underline{\operatorname{Aut}_\otimes^{\rm br}(\cB)}$}.

\subsection[The 2-categorical group of G-crossed braided extensions]{The 2-categorical group of $\boldsymbol{G}$-crossed braided extensions}

In a $G$-crossed braided fusion category $\cross{\cB}{G} = \bigoplus_{g \in G} \cB_g$, the trivial component $\cB_e$ is a braided fusion subcategory, each component $\cB_g$ is an invertible $\cB_e$-module category, and the functors~${M_{g,h}\colon \cB_g \boxtimes_{\cB_e} \cB_h \to \cB_{gh}}$ induced from the tensor product are $\cB_e$-module equivalences \cite[Theorem 6.1]{ENO3}.

The following theorem establishes a crucial connection between $G$-crossed braided categories and morphisms of 3-groups.

\begin{Theorem}[{\cite[Theorem 7.12]{ENO3}}] \label{theorem: levantamientos y g-crossed}
Let $\cB$ be a braided fusion category. Equivalence classes of braided $G$-crossed categories $\cross{\cB}{G}$ with a faithful $G$-grading and trivial component $\cB$ are in bijection with morphisms of $3$-groups $\underline{\underline{G}} \to \underline{\underline{\operatorname{Pic}(\cB)}}$.
\end{Theorem}

\subsection[The H\^4-obstruction]{The $\boldsymbol{H^4}$-obstruction}\label{subsec: H4 obstrucion}

We are primarily interested in fermionic actions \smash{$\widetilde{\rho}\colon \underline{G} \to \underline{\operatorname{Aut}_\otimes^{\rm br}(\cB, f)}$} that can be lifted to homomorphisms of 3-groups \smash{$\widetilde{\widetilde{\rho}}\colon \underline{\underline{G}} \to \underline{\underline{\operatorname{Pic}(\cB)}}$}. The obstruction to the existence of such liftings is measured by a cohomological invariant in $H^4\bigl(G, \mathbb{C}^\times\bigr)$.

For a non-degenerate braided fusion category $\cB$, by \cite[Section 6.1]{GELAKI20081053}, there is a canonical group isomorphism
\[
\Xi\colon\ \operatorname{Inv}(\cB) \to \widehat{K_0(\cB)}, \qquad \Xi(X_i)(X_j)=c_{X_j,X_i}\circ c_{X_i,X_j}.
\]

Let $\widetilde{\rho}$ be a categorical $G$-action determined by data $(g_*, \psi^g, \varphi_{g,h})\colon \underline{G} \to \underline{\operatorname{Aut}_\otimes^{\rm br}(\cB)}$. It follows from the isomorphism $\Xi$ and Proposition \ref{Proposition:Obstruction-bosonic} that the equivalence classes of such categorical liftings form a torsor over $H^2_\rho(G, \operatorname{Inv}(\cB))$. Given a 2-cocycle $\mu \in Z^2_\rho(G, \operatorname{Inv}(\cB))$, we denote by~${\mu \triangleright \widetilde{\rho}}$ the twisted categorical action.

Suppose that $\widetilde{\rho}$ admits a lifting $\widetilde{\widetilde{\rho}}$, that is, a homomorphism of $3$-groups \smash{$\widetilde{\widetilde{\rho}}\colon \underline{\underline{G}} \to \uu{\Pic{\cB}}$} whose truncation is $\widetilde{\rho}$. The question now is: given $\mu \in Z^2_\rho(G, \operatorname{Inv}(\cB))$ when can the new categorical action $\mu \triangleright \widetilde{\rho}$ be also lifted to a 3-group homomorphism? The answer is provided by the vanishing of a certain 4-cocycle, which we now define.

\begin{Definition}\label{def: formula H4}
The \emph{$H^4$-obstruction} of a pair $(\widetilde{\rho}, \mu)$ is the 4-cocycle $O_4(\widetilde{\rho}, \mu) \in H^4\bigl(G, \mathbb{C}^\times\bigr)$ described by
\begin{align}
O_4(\tilde{\rho}, \mu) ={} & \varphi_{g,h}(\mu(h,k)) \cdot c((gh)_*\mu(k,l), \mu(g,h)) \nonumber\\
& \times \frac{\psi^g(\mu(h,k), \mu(hk,l)) \cdot \omega((gh)_*\mu(k,l), \mu(g,h), \mu(gh,kl))}{(\psi^g)(h_*\mu(k,l), \mu(h,kl)) \cdot \omega((gh)_*\mu(k,l), g_*\mu(h,kl), \mu(g,hkl))} \nonumber \\
& \times \frac{\omega(\mu(g,h), \mu(gh,k), \mu(ghk,l)) \cdot \omega(g_*\mu(h,k), g_*\mu(hk,l), \mu(g,hkl))}{\omega(\mu(g,h), (gh)_*\mu(k,l), \mu(gh,kl)) \cdot \omega(g_*\mu(h,k), \mu(g,hk), \mu(ghk,l))},\label{equation: h4-obstruction}
\end{align}
where $\omega$ is the associativity constraint of the category and $c$ is the braiding.
\end{Definition}
\begin{Proposition}[{\cite[Proposition 9]{SCJZ}}]\label{prop: formula H4 anomaia}
If $\cB$ is a non-degenerate braided fusion category, the homomorphism of $2$-groups $(\mu \triangleright \widetilde{\rho})\colon \underline{G} \to \underline{\operatorname{Pic}(\cB)}$ can be lifted to a $3$-group homomorphism if and only if $O_4(\widetilde{\rho}, \mu)$ defined by \eqref{equation: h4-obstruction} is trivial in cohomology.
\end{Proposition}

\subsection{Non-degenerate extensions of braided fusion categories}\label{subsection: non-degenerate extension of braided fusion categories}

The concept of minimal modular extension of a braided fusion category was introduced by M\"uger in \cite{Mu2}. Here, we present a natural generalization where only the non-degeneracy condition is essential.

\begin{Definition}\label{def: minimal non-degenerate extension}
Let $(\cB,c)$ be a braided fusion category. A \emph{minimal non-degenerate extension} of $\cB$ is a pair $(\cM,i)$, where $\cM$ is a non-degenerate braided fusion category such that
\[\FPdim(\cM)=\FPdim(\cB)\FPdim(\mathcal{Z}_2(\cB))\]
and $i\colon \cB \to \cM$ is a braided full embedding.

Two minimal non-degenerate extensions $(\cM,i)$ and $(\cM',i')$ are equivalent if there exists a~braided equivalence $F\colon \cM \to \cM'$ such that $F \circ i \simeq i'$.
\end{Definition}

Every unitary braided fusion category admits a unique unitary ribbon structure \cite[Theorem~3.5]{Galindo-unitaria}. Therefore, minimal non-degenerate \emph{unitary} extensions of a unitary braided fusion category are modular extensions in the sense of \cite{Mu2}.

Weakly group-theoretical fusion categories admit a unique unitary structure (see \cite[Theo\-rem~5.20]{GHR}). Consequently, for weakly group-theoretical braided fusion categories, every minimal non-de\-gen\-erate extension is equivalent to a unitary minimal modular extension.

In \cite{LKW}, the set of equivalence classes of minimal non-degenerate extensions of a unitary braided fusion category $\cB$ is denoted by $\mathcal{M}_{\rm ext}(\cB)$. When $\cB$ is a symmetric fusion category, $\mathcal{M}_{\rm ext}(\cB)$ forms an abelian group.

\begin{Example}[minimal non-degenerate extensions of Tannakian fusion categories]\label{bosonicsymmetry}
For symmetric Tannakian categories $\Rep(G)$, the group of minimal non-degenerate extensions (up to equivalence) is isomorphic to $H^3(G,\mathbb{C}^\times)$. For each $\omega \in H^3(G,\mathbb{C}^\times)$, the Drinfeld center of the category of $G$-graded vector spaces with associativity constraint defined by $\omega$, denoted $\mathcal{Z}\bigl(\operatorname{Vec}_G^\omega\bigr)$, is a minimal non-degenerate extension of $\Rep(G)$ (see \cite{LKW}).
\end{Example}

\begin{Example}[minimal modular extensions of $\SV$]\label{example:modular-extensions-SV}
The symmetric super-Tannakian category $\SV$ has 16 minimal non-degenerate extensions (up to equivalence), which can be classified into two families:
\begin{enumerate}\itemsep=0pt
 \item[(1)] Eight Ising braided categories parametrized by $\zeta$, where $\zeta^8 = -1$. These correspond to the different central charges of the Ising braided categories.
 \item[(2)] Eight pointed non-degenerate spin-braided fusion categories $\operatorname{Vec}_A^q$ where $A$ is an abelian group of order four and $q\colon A \to \mathbb{Q}/\mathbb{Z}$ is a non-degenerate quadratic form, see Example~\ref{ex: pointed 4 rank}.
\end{enumerate}

For more details, see \cite{DGNO,Kitaev20062}.
\end{Example}

\subsection{Obstruction theory for minimal non-degenerate extensions}\label{subsection: Obstruction theory to existence of minimal non-degenerate extensions}

Following \cite{Bruguières2000,Galois-Mueger}, we call a braided fusion category \emph{modularizable} if $\mathcal{Z}_2(\cB)$ is Tannakian.

\begin{Definition}\label{def: H4 anolmaly}
Let $\cB$ be a modularizable braided fusion category with $\mathcal{Z}_2(\cB) = \Rep(G)$. By \cite[Proposition 4.30\,(iii)]{DGNO}, the de-equivariantization $\cB_G$ is non-degenerate. The category $\cB_G$ carries a canonical action of $G$, which induces a monoidal functor
\[
\xymatrix{
 \underline{G} \ar[r]^-{\widetilde{\rho}} & \underline{\operatorname{Aut}_\otimes^{\rm br}(\cB_{G})} \ar[r]^{\Phi} & \underline{\operatorname{Pic}(\cB_G)},
}
\]
where $\Phi$ is the equivalence from \cite[Theorem 5.2]{ENO3}. We define the \emph{$H^4$-anomaly} of $\cB$, denoted~$O_4(\cB)$, as the $H^4$-obstruction of the composition $\Phi \circ \widetilde{\rho}$ in $H^4(G,\mathbb{C}^\times)$.
\end{Definition}

\begin{Definition}[{\cite{ENO2}}]
A braided fusion category $\cB$ is called \emph{slightly degenerate} if $\mathcal{Z}_2(\cB)$ is braided equivalent to $\SV$.
\end{Definition}

If $\cB$ is non-modularizable, that is, $\mathcal{Z}_2(\cB) \cong \Rep\bigl(\widetilde{G},z\bigr)$, then the maximal central Tannakian subcategory of $\cB$ is equivalent to $\Rep(G)$ with \smash{$G \cong \widetilde{G}/\langle z \rangle$}.

The following theorem provides necessary and sufficient conditions for the existence of minimal non-degenerate extensions.

\begin{Theorem}\label{theorem:teorema de-equivariantización y extensiones}
Let $\cB$ be a braided fusion category with non-trivial maximal central Tannakian subcategory $\Rep(G) \subseteq \mathcal{Z}_2(\cB)$.
\begin{enumerate}\itemsep=0pt
 \item[$(1)$]  If $\cB$ is modularizable, then $\cB$ admits a minimal non-degenerate extension if and only if the $H^4$-anomaly $O_4(\cB)$ vanishes.

 \item[$(2)$] If $\cB$ is non-modularizable with $\mathcal{Z}_2(\cB) = \Rep\bigl(\widetilde{G},z\bigr)$, then $\cB$ admits a minimal non-de\-gen\-erate extension if and only if the following conditions hold:
 \begin{enumerate}\itemsep=0pt
 \item[$(a)$] The slightly degenerate braided fusion category $\cB_G$ has a minimal non-degenerate extension $\cS$.

 \item[$(b)$] There exists a fermionic action of $\bigl(\widetilde{G},z\bigr)$ on $\cS$ such that $\cB_G$ is $G$-stable, and the restriction to $\cB_G$ coincides with the canonical action of $G$ on $\cB_G$.

 \item[$(c)$] The anomaly $O_4\bigl(\cS^G\bigr)$ vanishes.
 \end{enumerate}
\end{enumerate}
\end{Theorem}

\begin{proof}
We prove each part separately using the relationship between centralizers and (de)equiv\-ari\-antization established in Section~\ref{subsection: centralizer}.

(1) Suppose first that $\cB$ admits a minimal non-degenerate extension $\cM$. Since $\mathcal{Z}_2(\cB) = \Rep(G)$, we have $\cB \subset \cM$ and $\Rep(G) \subseteq C_\cM(\Rep(G))$.

By \cite[Theorem 4.4]{DGNO}, there is a canonical braided inclusion $\cB_G \subseteq (C_\cM(\Rep(G)))_G$. By~\cite[Proposition 4.56\,(i)]{DGNO}, we have $(C_\cM(\Rep(G)))_G = (\cM_G)_e$, where $(\cM_G)_e$ denotes the trivial component of the braided $G$-crossed category $\cM_G$.

We claim that $\cB_G = (\cM_G)_e$. To see that, let us compute the Frobenius--Perron dimensions. Since $\cM$ is a minimal non-degenerate extension of $\cB$, we have
\begin{align*}
\FPdim((\cM_G)_e) &= \frac{\FPdim(\cM)}{|G|^2} = \frac{\FPdim(\cB) \cdot |\mathcal{Z}_2(\cB)|}{|G|^2}\\ &= \frac{\FPdim(\cB) \cdot |G|}{|G|^2} = \frac{\FPdim(\cB)}{|G|} = \FPdim(\cB_G).
\end{align*}

Since $\cB_G \subseteq (\cM_G)_e$ and both categories have the same Frobenius--Perron dimension, we conclude $\cB_G = (\cM_G)_e$. This means that $\cM_G$ is a braided $G$-crossed extension of $\cB_G$ with $\cB = (\cB_G)^G \subset \cM$. By Theorem~\ref{theorem: levantamientos y g-crossed} and Proposition \ref{prop: formula H4 anomaia}, the existence of such a $G$-crossed extension is equivalent to the vanishing of the $H^4$-anomaly $O_4(\cB)$.

Conversely, assume $O_4(\cB) = 0$. By Theorem~\ref{theorem: levantamientos y g-crossed}, there exists a braided $G$-crossed extension~$\cL$ of $\cB_G$. Taking the $G$-equivariantization, we obtain $\cB = (\cB_G)^G \subset \cL^G$.

By \cite[Proposition 4.56\,(ii)]{DGNO}, $\cL^G$ is non-degenerate since $\cB_G$ is non-degenerate and the $G$-grading on $\cL$ is faithful. Moreover,
\begin{align*}
\FPdim\bigl(\cL^G\bigr) &= |G| \cdot \FPdim(\cL) = |G|^2 \cdot \FPdim(\cB_G) \\&= |G| \cdot \FPdim(\cB) = \FPdim(\cB) \cdot |\mathcal{Z}_2(\cB)|.
\end{align*}

Therefore, $\cL^G$ is a minimal non-degenerate extension of $\cB$.

(2) Let $\cB$ be a braided fusion category with $\mathcal{Z}_2(\cB) = \Rep\bigl(\widetilde{G},z\bigr)$ and $G = \widetilde{G}/\langle z \rangle$.

Suppose first that $\cB$ admits a minimal non-degenerate extension $\cM$. By \cite[Proposition~4.30 (iii)]{DGNO}, the de-equivariantization $\cB_G$ is slightly degenerate with $\mathcal{Z}_2(\cB_G) \cong \SV$. We have $\cB_G \subseteq (\cM_G)_e$ and $(\cM_G)_e$ is non-degenerate by \cite[Proposition 4.56\,(ii)]{DGNO}.

Let us verify that $\cS := (\cM_G)_e$ is a minimal non-degenerate extension of $\cB_G$. Since $\cM$ is a~minimal extension of $\cB$, we have
\begin{align*}
 \FPdim((\cM_G)_e) &= \frac{\FPdim(\cM)}{|G|^2} = \frac{\FPdim(\cB) \cdot |\mathcal{Z}_2(\cB)|}{|G|^2}= \frac{\FPdim(\cB) \cdot 2|G|}{|G|^2}\\
 & = \frac{2\FPdim(\cB)}{|G|} = 2 \cdot \FPdim(\cB_G) = \FPdim(\cB_G) \cdot |\mathcal{Z}_2(\cB_G)|.
\end{align*}

This establishes condition (2a). Since $\cB_G$ is slightly degenerate, $\cS$ has a canonical spin-braided structure with fermion given by the generator of $\mathcal{Z}_2(\cB_G) \cong \SV$.

For condition~(2b), we observe that $\Rep\bigl(\widetilde{G},z\bigr) \subset \cB \subset \cS^G$. By Corollary~\ref{coro: equivalencia dequivariantizacion- equivariantizacion braided case}, this implies that \smash{$\bigl(\widetilde{G},z\bigr)$} acts by braided autoequivalences on $\cS$, the subcategory $\cB_G$ is $G$-stable under this action, and the restriction coincides with the canonical action of $G$ on $\cB_G$.

For condition (2c), since $\cM$ is a minimal non-degenerate extension of $\cS^G$, the vanishing of~$O_4\bigl(\cS^G\bigr)$ follows from part~(1).

Conversely, assume conditions (2a), (2b), and (2c) hold. By condition (a), we have a~minimal non-degenerate extension $\cS$ of $\cB_G$. By condition (2b), there is a canonical braided inclusion $\cB_G \subset \cS$, which by equivariantization gives $\cB = (\cB_G)^G \subset \cS^G$. By condition~(2c) and part~(1), $\cS^G$ has a minimal non-degenerate extension $\cM$.

It suffices to check that $\cM$ is a minimal non-degenerate extension of $\cB$. Since $\cS^G$ is modularizable and $\cM$ is its minimal extension, we have
\begin{align*}
\FPdim(\cM) &= \FPdim\bigl(\cS^G\bigr) \cdot |\mathcal{Z}_2\bigl(\cS^G\bigr)| = |G|^2 \cdot \FPdim(\cS) \cdot |G|= |G|^2 \cdot 2 \cdot \FPdim(\cB_G) \\
&= 2|G|^2 \cdot \frac{\FPdim(\cB)}{|G|} = 2|G| \cdot \FPdim(\cB) = \FPdim(\cB) \cdot |\mathcal{Z}_2(\cB)|.
\end{align*}

Therefore, $\cM$ is a minimal non-degenerate extension of $\cB$.
\end{proof}

If $\Rep(G)$ is a central Tannakian subcategory of a braided fusion category $\cB$, then by the proof of Theorem~\ref{theorem:teorema de-equivariantización y extensiones}, the de-equivariantization defines a map
\begin{align*}
D\colon \ \mathcal{M}_{\text{ext}}(\cB) \to \mathcal{M}_{\text{ext}}(\cB_G) ,\qquad
\cM \mapsto (\cM_G)_e.
\end{align*}
This map was also defined in \cite[Section~5.2]{LKW}, and it was shown that when $\cB$ is symmetric, it is a homomorphism of abelian groups. In particular, when $\cB = \Rep\bigl(\widetilde{G}, z\bigr)$, we obtain a group homomorphism
\begin{align}\label{eq: map D}
D\colon \ \mathcal{M}_{\text{ext}}\bigl(\Rep\bigl(\widetilde{G}, z\bigr)\bigr) \to \mathcal{M}_{\text{ext}}(\SV).
\end{align}

The following corollary characterizes when this map is surjective.

\begin{Corollary}\label{corol:Ostrik}
Let $\Rep\bigl(\widetilde{G},z\bigr)$ be a super-group. The map \[
D\colon\ \mathcal{M}_{\rm ext}\bigl(\Rep\bigl(\widetilde{G},z\bigr)\bigr) \to \mathcal{M}_{\rm ext}(\SV)\]
 is surjective if and only if \smash{$\bigl(\widetilde{G},z\bigr)$} is a trivial super-group.
\end{Corollary}

\begin{proof}
If $D$ is surjective, then by Example \ref{example:modular-extensions-SV}, there exists $\cM \in \mathcal{M}_{\rm ext}\bigl(\Rep\bigl(\widetilde{G},z\bigr)\bigr)$ such that~$(\cM_G)_e$ is a spin-braided Ising category and \smash{$\bigl(\widetilde{G},z\bigr)$} acts on $(\cM_G)_e$. By Proposition \ref{proposition: fermionicactionsonIsing}, we conclude that \smash{$\bigl(\widetilde{G},z\bigr)$} is a trivial super-group.

Conversely, if $\widetilde{G} = G \times \mathbb{Z}/2\mathbb{Z}$ with central involution $z = (e,1)$, then $\Rep\bigl(\widetilde{G},z\bigr) = \Rep(G) \boxtimes \SV$ as braided fusion categories. For any $\mathcal{M} \in \mathcal{M}_{\rm ext}(\SV)$, the non-degenerate fusion category~${\mathcal{Z}(\Rep(G)) \boxtimes \cM \supset \Rep(G) \boxtimes \SV}$ is a minimal non-degenerate extension with
\begin{equation*}
D(\mathcal{Z}(\Rep(G)) \boxtimes \cM) = ((\mathcal{Z}(\Rep(G)) \boxtimes \cM)_G)_e = \cM.
\end{equation*}
Therefore, the map $D$ is surjective.
\end{proof}

\subsection[Examples of modularizable braided fusion categories without minimal extensions]{Examples of modularizable braided fusion categories\\ without minimal extensions}\label{subsection:obstruction 4}

Let $G$ be a finite group and $\cB$ a non-degenerate braided fusion category. If we consider the trivial 3-homomorphism \smash{$\widetilde{\widetilde{\text{triv}}}\colon\underline{\underline{G}} \to \underline{\underline{\operatorname{Pic}(\cB)}}$}, then the associated categorical action is also trivial. For~${\mu \in Z^2(G, \operatorname{Inv}(\cB))}$, we ask when the braided fusion category $\cB^G$ obtained by equivariantization associated to the categorical action $\mu \triangleright \widetilde{\text{triv}}$ admits a minimal non-degenerate extension.

\begin{Corollary}\label{formula 4 quasi-trivial}
Let $G$ be a finite group and $(\cB, \otimes, \mathbf{1}, \omega, c)$ be a non-degenerate braided fusion category. The equivariantization $\cB^G$ associated with $\mu \in H^2(G, \operatorname{Inv}(\cB))$ has a minimal non-degenerate extension if and only if the cohomology class of the $4$-cocycle 	
\begin{align}
O_4(\tilde{\rho}, \mu) ={} & c(\mu(k,l), \mu(g,h)) \cdot \frac{\omega(\mu(k,l), \mu(g,h), \mu(gh,kl))}{\omega(\mu(k,l), \mu(h,kl), \mu(g,hkl))}\nonumber \\
& \cdot \frac{\omega(\mu(g,h), \mu(gh,k), \mu(ghk,l)) \cdot \omega(\mu(h,k), \mu(hk,l), \mu(g,hkl))}{\omega(\mu(g,h), \mu(k,l), \mu(gh,kl)) \cdot \omega(\mu(h,k), \mu(g,hk), \mu(ghk,l))}\label{formula o4}
\end{align}
vanishes in $H^4\bigl(G, \mathbb{C}^\times\bigr)$.
\end{Corollary}

\begin{proof}
Since $\cB$ is non-degenerate, $\cD := \cB^G$ is modularizable. By Theorem~\ref{theorem:teorema de-equivariantización y extensiones}\,(1), $\cD$ admits a~minimal non-degenerate extension if and only if its $H^4$-anomaly vanishes. By Proposition~\ref{prop: formula H4 anomaia}, the $H^4$-anomaly of $\cD$ corresponds to the obstruction $O_4\bigl(\widetilde{\operatorname{triv}}, \mu\bigr)$ for lifting the categorical action~\smash{$\mu \triangleright \widetilde{\operatorname{triv}}$}, which is precisely the formula given in~\eqref{formula o4}.
\end{proof}

Using Corollary~\ref{formula 4 quasi-trivial}, we can construct explicit examples of modularizable braided fusion categories without minimal non-degenerate extensions. Our approach relies on computing 4-cocycles whose cohomology classes are non-trivial.

For a direct product of groups $G = A \times B$, we recall that there is a canonical filtration in group cohomology
\[H^n = \mathfrak{F}_0(H^n) \supset \mathfrak{F}_1(H^n) \supset \cdots \supset \mathfrak{F}_{n-1}(H^n) \supset \mathfrak{F}_n(H^n) \supset 0, \]
where $H^n = H^n(A \times B, M)$. This filtration gives rise to group homomorphisms
\begin{align*}
P_k\colon\ \operatorname{gr}_k(H^n(A \times B, M)) \to H^k\bigl(A, H^{n-k}(B, M)\bigr).
\end{align*}
We will use these maps to detect non-trivial cohomology classes in $H^4\bigl(G, \mathbb{C}^\times\bigr)$.

The following example was first discovered by Drinfeld in unpublished notes. We identify the group of roots of unity in $\mathbb{C}$ with $\mathbb{Q}/\mathbb{Z}$ throughout.

\begin{Proposition}
Let $(\omega, c) \in Z^3_{ab}(\mathbb{Z}/2\mathbb{Z}, \mathbb{Q}/\mathbb{Z})$ be the non-degenerate abelian $3$-cocycle given~by
\begin{align*}
 c(x,y) = \frac{xy}{4}, \qquad \omega(x,y,z) = \frac{xyz}{2},\qquad x,y,z\in\{0,1\}
\end{align*}
and $\mu \in Z^2(\mathbb{Z}/2\mathbb{Z} \times \mathbb{Z}/2\mathbb{Z}, \mathbb{Z}/2\mathbb{Z})$ given by $\mu((a_1,b_1),(a_2,b_2)) = a_1b_2$. Then the braided fusion category \smash{$\bigl(\operatorname{Vec}_{\mathbb{Z}/2\mathbb{Z}}^{(\omega,c)}\bigr)^{\mathbb{Z}/2\mathbb{Z} \times \mathbb{Z}/2\mathbb{Z}}$} associated to the categorical action $\mu \triangleright \widetilde{\operatorname{triv}}$ does not admit a minimal non-degenerate extension.
\end{Proposition}

\begin{proof}

Let $G = \mathbb{Z}/2\mathbb{Z} \times \mathbb{Z}/2\mathbb{Z}$. We adopt the convention that every element $g_i \in G$ will be written as $g_i = (a_i, b_i)$.

After simplification, the 4-cocycle $O_4$ defined by equation \eqref{formula o4} is given by
\[
O_4(g_1,g_2,g_3,g_4) = \frac{a_1b_2a_3b_4}{4} + \frac{a_1a_2a_3b_2b_4 + a_1a_2a_3b_4 + a_1b_2b_3b_4 + a_1a_3b_2b_3b_4}{2}.
\]

It follows from Corollary~\ref{formula 4 quasi-trivial} that if the cohomology class of $O_4$ is non-zero, the braided fusion category \smash{$\bigl(\operatorname{Vec}_{\mathbb{Z}/2\mathbb{Z}}^{(\omega,c)}\bigr)^{\mathbb{Z}/2\mathbb{Z} \times \mathbb{Z}/2\mathbb{Z}}$} does not admit a minimal non-degenerate extension.

We define the 3-cochain
$
 p(g_1,g_2,g_3) = \frac{a_1a_2b_3}{8} - \frac{a_1b_2a_2b_3}{4}$,
and compute
\begin{align*}
 (\partial(p) - O_4)(g_1,g_2,g_3,g_4) = \frac{a_1a_2a_3b_4 + a_1b_2b_3b_4}{2}.
\end{align*}

Therefore, $O_4$ is cohomologous to
$
 \widetilde{O_4}(g_1,g_2,g_3,g_4) = \frac{a_1a_2a_3b_4 + a_1b_2b_3b_4}{2}$.

Thus \smash{$\big[\widetilde{O_4}\big] \in \mathfrak{F}_1\bigl(H^4(A \times B, \mathbb{Q}/\mathbb{Z})\bigr)$} and
\[P_1\bigl(\big[\widetilde{O_4}\big]\bigr) \in \operatorname{Hom}\bigl(\mathbb{Z}/2\mathbb{Z}, H^3(\mathbb{Z}/2\mathbb{Z}, \mathbb{Q}/\mathbb{Z})\bigr) \cong \mathbb{Z}/2\mathbb{Z}\]
is induced by the map
\begin{align*}
 \mathbb{Z}/2\mathbb{Z} \to Z^3(\mathbb{Z}/2\mathbb{Z}, \mathbb{Q}/\mathbb{Z}), \qquad a_1 \mapsto \left[(b_2,b_3,b_4) \mapsto \frac{a_1b_2b_3b_4}{2}\right].
\end{align*}

Since the cohomology class of the 3-cocycle $\alpha(b_2,b_3,b_4) = \frac{b_2b_3b_4}{2}$ is non-trivial, we have $P_1\bigl(\big[\widetilde{O_4}\big]\bigr) \neq 0$. Hence, $[O_4] \neq 0$ in $H^4(G, \mathbb{Q}/\mathbb{Z})$, which proves the result.
\end{proof}

\begin{Proposition}
Let $m$ be an odd integer, $(0,c) \in Z^3_{ab}(\mathbb{Z}/m\mathbb{Z}, \mathbb{Q}/\mathbb{Z})$ be the non-degenerate abelian $3$-cocycle given by
$c(x,y) = \frac{xy}{m}$,
and $\mu \in Z^2\bigl((\mathbb{Z}/m\mathbb{Z})^{\oplus n} \times (\mathbb{Z}/m\mathbb{Z})^{\oplus n}, \mathbb{Z}/m\mathbb{Z}\bigr)$ given by
\[\mu\bigl(\bigl(\va_1,\vb_1\bigr),\bigl(\va_2,\vb_2\bigr)\bigr) = \sum_{i=1}^{n} a_1^i b_2^i,\]
with $n \geq 2$. Then the braided fusion category \smash{$\bigl(\operatorname{Vec}_{\mathbb{Z}/m\mathbb{Z}}^{c}\bigr)^{(\mathbb{Z}/m\mathbb{Z})^{\oplus n} \times (\mathbb{Z}/m\mathbb{Z})^{\oplus n}}$}, associated to the categorical action \smash{$\mu \triangleright \widetilde{\operatorname{triv}}$} does not admit a minimal non-degenerate extension..
\end{Proposition}

\begin{proof}
Let $A = B = (\mathbb{Z}/m\mathbb{Z})^{\oplus n}$ and $G = A \times B$. We adopt the convention that every element~${\vg_i \in G}$ will be written as $\vec{g}_i = (\va_i, \vb_i)$, where $\va_i \in A$ and $\vb_i \in B$.

Given $\va \in A$ and $\vb \in B$, we define $\va\vb = \sum_{i=1}^n a_i b_i \in \mathbb{Z}/m\mathbb{Z}$. Then $\mu \in Z^2(G, \mathbb{Z}/m\mathbb{Z})$ is defined by \smash{$\mu(\vg_1, \vg_2) = \va_1\vb_2$}.
The 4-cocycle $O_4$ defined by equation \eqref{formula o4} is given by
\[
 O_4(\vg_1,\vg_2,\vg_3,\vg_4) = \frac{\bigl(\va_1\vb_2\bigr)\bigl(\va_3\vb_4\bigr)}{m}.
 \]

It follows from Corollary~\ref{formula 4 quasi-trivial} that if the cohomology class of $O_4$ is non-zero, the braided fusion category \smash{$\bigl(\operatorname{Vec}_{\mathbb{Z}/m\mathbb{Z}}^{c}\bigr)^{(\mathbb{Z}/m\mathbb{Z})^{\oplus n} \times (\mathbb{Z}/m\mathbb{Z})^{\oplus n}}$} does not admit a minimal non-degenerate extension.

We define the 3-cochain
\smash{$p(\vg_1,\vg_2,\vg_3) = \frac{(\va_1\vb_2)(\va_2\vb_3)}{m}$}
such that
\begin{align*}
 \partial(p)(\vg_1,\vg_2,\vg_3,\vg_4) = \frac{\bigl(\va_1\vb_2\bigr)\bigl(\va_3\vb_4\bigr)}{m} + \frac{\bigl(\va_1\vb_3\bigr)\bigl(\va_2\vb_4\bigr)}{m}.
\end{align*}
Therefore, $O_4$ is cohomologous to
\smash{$\widetilde{O}_4(\vg_1,\vg_2,\vg_3,\vg_4) \!=\! -\frac{(\va_1\vb_3)(\va_2\vb_4)}{m}$}.
Hence \smash{$\big[\widetilde{O}_4\big] \!\in\! \mathfrak{F}_2\bigl(H^4(G, \mathbb{Q}/\mathbb{Z})\bigr)$} and
\smash{$P_2\bigl(\big[\widetilde{O_4}\big]\bigr) \in H^2\bigl(A, H^2(B, \mathbb{Q}/\mathbb{Z})\bigr)$}
is induced by the 2-cocycle $\psi \in Z^2\bigl(A, Z^2(B, \mathbb{Z}/m\mathbb{Z})\bigr)$, where
\[\psi(\va_1,\va_2)\bigl(\vb_1,\vb_2\bigr) = \widetilde{O}_4\bigl(\va_1,\va_2,\vb_1,\vb_2\bigr) = -\frac{\bigl(\va_1\vb_1\bigr)\bigl(\va_2\vb_2\bigr)}{m}.\]

For any abelian group $\mathcal{A}$ and trivial $\mathcal{A}$-module $M$, the map
\begin{align*}
Z^2(\mathcal{A}, M) &\to \operatorname{Hom}\bigl(\wedge^2\mathcal{A}, M\bigr),\qquad
\alpha \mapsto [a_1 \wedge a_2 \mapsto \alpha(a_1, a_2) - \alpha(a_2, a_1)]
\end{align*}
induces a group homomorphism
$
\operatorname{Alt}_{\mathcal{A}}\colon H^2(\mathcal{A}, M) \to \operatorname{Hom}\bigl(\wedge^2\mathcal{A}, M\bigr)$.

The map
\smash{$\widehat{\psi} := \operatorname{Alt}_{A}(\operatorname{Alt}_{B}(\psi)) \in \operatorname{Hom}\bigl(\wedge^2 A, \operatorname{Hom}(\wedge^2 B, \mathbb{Q}/\mathbb{Z})\bigr)$}
given by
\[\widehat{\psi}(\va_1,\va_2)\bigl(\vb_1,\vb_2\bigr) = 2\frac{\bigl(\va_1\vb_1\bigr)\bigl(\va_2\vb_2\bigr) - \bigl(\va_1\vb_2\bigr)\bigl(\va_2\vb_1\bigr)}{m}\]
depends only on the cohomology class of $O_4$.

Let $\vec{e}_1=(1,0,\ldots,0)$ and $\vec{e}_2=(0,1,0,\ldots,0)$. Since $\widehat{\psi}(\vec{e}_1, \vec{e}_2)(\vec{e}_1, \vec{e}_2) = \frac{2}{m} \neq 0$, the map $\widehat{\psi}$ is non-trivial, and then the cohomology class of $O_4$ is also non-trivial.
\end{proof}

\subsection[H\^4-obstruction for non-degenerate braided categories of rank four]{$\boldsymbol{H^4}$-obstruction for non-degenerate braided categories of rank four}

We have studied the fermionic actions of super-groups $(G, \alpha)$ and the $\alpha$-liftings of homomorphisms $\rho\colon G \to \operatorname{Aut}_\otimes^{\rm br}(\cB, f)$ for non-degenerate spin-braided fusion categories of rank four. In this section, we present formulas for the $H^4$-obstruction to the existence of a 3-group homomorphism
\[\widetilde{\widetilde{\rho}}\colon\ \underline{\underline{G}} \to \underline{\underline{\operatorname{Aut}_\otimes^{\rm br}(\cB)}} \cong \underline{\underline{\operatorname{Pic}(\cB)}}\]
when $\cB$ is a non-degenerate spin-braided fusion category of rank four (see Examples \ref{ex: pointed 4 rank} and~\ref{Example:Ising}).

This topic is important for the construction of minimal non-degenerate extensions of super-groups. We recall that there is a canonical map
\[D\colon\ \mathcal{M}_{\text{ext}}\bigl(\Rep\bigl(\widetilde{G}, z\bigr)\bigr) \to \mathcal{M}_{\text{ext}}(\SV)\]
defined in \eqref{eq: map D}.

The problem of determining which minimal non-degenerate extensions of $\SV$ lie in the image of $D$ reduces to the following computational steps:
\begin{enumerate}\itemsep=0pt
\item[(1)] Classify all fermionic actions of $\bigl(\widetilde{G},z\bigr)$ on pointed non-degenerate spin-braided categories of rank four (see Theorem~\ref{theorem: extensions of A with order 4} for the pointed and Proposition \ref{proposition: fermionicactionsonIsing} for the Ising case).
\item[(2)] For each such action, compute the associated $H^4$-obstruction.
\item[(3)] Verify whether this obstruction vanishes in cohomology.
\end{enumerate}
If the $H^4$-obstruction vanishes, then the corresponding non-degenerate spin-braided fusion category lies in the image of $D$, and we can construct the corresponding minimal non-degenerate extensions of $\Rep\bigl(\widetilde{G},z\bigr)$ by taking $G$-equivariantization of the associated braided $G$-crossed fusion category.

To apply the $H^4$-obstruction formula from Definition~\ref{def: formula H4} to any lifting of a group homomorphism $\xi\colon G \to \operatorname{Aut}_\otimes^{\rm br}(\cB, f)$, we need to construct a 3-group homomorphism \smash{$\widetilde{\widetilde{\xi}}\colon \underline{\underline{G}} \to \underline{\underline{\operatorname{Pic}(\cB)}}$}. This construction was straightforward for trivial actions, allowing us to perform explicit computations. The following result shows that when $\cB$ is a non-degenerate spin-braided fusion category of rank four, such lifting always exist.

\begin{Lemma}\label{extensiontriviasupergroup}
Let $(\cB, f)$ be a non-degenerate spin-braided fusion category of rank four and $G$ a finite group with a group homomorphism $\xi\colon G \to \operatorname{Aut}_\otimes^{\rm br}(\cB, f)$. Then $\xi$ can be extended to a~categorical action and to a $3$-group homomorphism \smash{$\widetilde{\widetilde{\xi}}\colon \underline{\underline{G}} \to \underline{\underline{\operatorname{Pic}(\cB)}}$}.
\end{Lemma}

\begin{proof}
If $\xi$ is trivial, we can always take the trivial lifting regardless of $G$ or $\cB$. Since all tensor automorphisms of the Ising fusion category are trivial, $\xi$ is trivial and we have the liftings in this case.

Consider the case where $\xi$ is non-trivial. Then $\cB$ is a pointed non-degenerate spin-braided fusion category. By Theorem~\ref{theorem: extensions of A with order 4}, we have $\operatorname{Aut}_\otimes^{\rm br}(\cB, f) \cong \mathbb{Z}/2\mathbb{Z}$, so $\xi$ is surjective and factors through an isomorphism $\rho\colon\mathbb{Z}/2\mathbb{Z} \to \operatorname{Aut}_\otimes^{\rm br}(\cB, f)$
\begin{equation*}
\xymatrix{
 G \ar[r]^-{\xi} \ar[d]_{\pi} & \operatorname{Aut}_\otimes^{\rm br}(\cB, f) \\
 \mathbb{Z}/2\mathbb{Z} \ar[ur]_{\rho}.
}
\end{equation*}

Our strategy is to first construct a 3-group lifting \smash{$\widetilde{\widetilde{\rho}}\colon \underline{\underline{\mathbb{Z}/2\mathbb{Z}}} \to \underline{\underline{\operatorname{Pic}(\cB)}}$} and then obtain the desired lifting of $\xi$ by composition with the projection $\pi$.

By Theorem~\ref{theorem: extensions of A with order 4}, the isomorphism $\rho$ always lifts to a categorical action $\widetilde{\rho}$. Since $H^4(\mathbb{Z}/2\mathbb{Z}, \allowbreak \mathbb{C}^\times)\allowbreak = 0$, the $H^4$-obstruction vanishes, so $\widetilde{\rho}$ lifts to a 3-group homomorphism \smash{$\widetilde{\widetilde{\rho}}$}.

To complete the construction, we define $\widetilde{\xi} = \widetilde{\rho} \circ \pi$ and \smash{$\widetilde{\widetilde{\xi}} = \widetilde{\widetilde{\rho}} \circ \pi$}. These compositions provide the required categorical action and 3-group homomorphism liftings of $\xi$.
\end{proof}

Using the lifting constructed in Lemma \ref{extensiontriviasupergroup}, we can now find a formula for the $H^4$-ob\-struc\-tion.

\begin{Theorem}
Let $(G, \alpha)$ be a super-group and $(\cB, f)$ a non-degenerate spin-braided fusion category of rank four. Let \smash{$\widetilde{\rho}_\alpha\colon \underline{G} \to \underline{\operatorname{Aut}_\otimes^{\rm br}(\cB, f)}$} be an $\alpha$-lifting with associated group homomorphism \smash{$\xi\colon G \to \operatorname{Aut}_\otimes^{\rm br}(\cB, f)$}. If \smash{$\widetilde{\widetilde{\xi}}\colon \underline{\underline{G}} \to \underline{\underline{\operatorname{Pic}(\cB)}}$} is a $3$-group lifting constructed by Lemma {\rm\ref{extensiontriviasupergroup}} and~${\mu \triangleright \widetilde{\xi} = \widetilde{\rho}_\alpha}$ for some \smash{$\mu \in Z^2_{\xi}(G, \operatorname{Inv}(\cB))$}, then $\widetilde{\rho}_\alpha$ can be extended to a $3$-group homomorphism if and only if the $4$-cocycle $O_4\bigl(\widetilde{\xi}, \mu\bigr)$ defined by \eqref{equation: h4-obstruction} is trivial in cohomology.
\end{Theorem}

\begin{proof}
In general, as we saw, if $\mu \in Z^2_{\xi}(G, \operatorname{Inv}(\cB))$, then $\mu$ defines a new categorical action~\smash{$\mu \triangleright \widetilde{\xi}$}. According to Proposition \ref{prop: formula H4 anomaia}, this categorical action has a lifting to a 3-homomorphism if and only if the 4-cocycle $O_4\bigl(\widetilde{\xi}, \mu\bigr)$ defined by \eqref{equation: h4-obstruction} is cohomologically trivial.
\end{proof}

\subsection*{Acknowledgements}

We are deeply grateful to the three anonymous referees for their thorough and insightful reports. Their numerous detailed comments and suggestions have significantly improved the clarity of exposition and overall quality of this work.

\pdfbookmark[1]{References}{ref}
\LastPageEnding

\end{document}